\documentclass[a4paper,10pt]{article}
\usepackage{amsmath,amsbsy,amsfonts,amssymb,amsthm,
amscd,amsxtra,amstext,amsopn,dsfont,yfonts}
\usepackage[english]{babel}
\usepackage{bbm,bm,makeidx,psfrag}
\usepackage{indentfirst,fancyhdr,float}
\usepackage{graphicx,epsfig}
\usepackage{caption}
\usepackage{latexsym}
\usepackage{amsfonts}
\usepackage{amsmath}
\usepackage{mathrsfs}
\usepackage[all]{xy}

\def\bes{\begin{equation*}}
\def\ees{\end{equation*}}
\def\ba{\begin{aligned}}
\def\ea{\end{aligned}}
\def\be{\begin{equation}}
\def\ee{\end{equation}}
\def\bc{\begin{cases}}
\def\ec{\end{cases}}

\newcommand{\zerarcounters}{\setcounter{equation}{0}}

\theoremstyle{plain}
\newtheorem{teo}{{\bf Theorem}}
\newtheorem{prop}{{\bf Proposition}}
\newtheorem{lemma}{{\bf Lemma}}

\newtheorem{hyp}{Hypothesis}

\theoremstyle{definition}

\newcommand{\prova}{\noindent{\it Proof. }}

\let\a=\alpha \let\b=\beta      \let\d=\delta     
\let\e=\varepsilon
  \let\h=\eta     \let\th=\theta \let\k=\kappa     

\let\m=\mu    \let\n=\nu           \let\p=\pi        

\let\s=\sigma \let\t=\tau         

   \let\o=\omega   

\let\G=\Gamma \let\D=\Delta   \let\Th=\Theta \let\L=\Lambda    

\let\P=\Pi      \let\F=\Phi

\def\PP{{\cal P}}\def\MM{{\cal M}} 
\def\CC{{\cal C}}\def\FF{{\cal F}}
\def\TT{{\cal T}}\def\NN{{\cal N}}
\def\RR{{\cal R}}
\def\GG{{\cal G}}
\def\QQ{{\cal Q}}

\def\aaa{{\mathfrak a}}

\def\eee{{\mathfrak e}}
\def\hhh{{\mathfrak h}}

\def\nnn{{\mathfrak n}}
\def\ppp{{\mathfrak p}}
\def\qqq{{\mathfrak q}}\def\rrr{{\mathfrak r}}
\def\sss{{\mathfrak s}}
\def\vvv{{\mathfrak v}}
\def\www{{\mathfrak w}}

\def\AAA{{\mathfrak A}}

\def\ie{{\emph{i.e.}\;}}

\def\hknu#1#2#3{{#1_{#3}^{[#2]}}}
\def\eknu#1#2#3{{#1_{#3}^{(#2)}}}

\def\dpr{\partial}

\def\Val{{\rm Val}}	\def\sign{{\rm sign}\,}
\def\deg{{\rm deg}}     \def\Syl{{\rm Syl}}

\def\ol#1{{\overline #1}}
\def\tilde#1{\widetilde{#1}}

\def\to{\rightarrow}

\def\la{\left\langle}
\def\ra{\right\rangle}
\def\mean#1{{\la #1\ra}}
\def\io{{\infty}}
\def\salto{{\vskip.2truecm}}
\def\bt{{\tilde{\b}}}

\def\RRR{{\mathbb R}}
\def\QQQ{{\mathbb Q}}
\def\ZZZ{{\mathbb Z}} 
\def\CCC{{\mathbb C}}
\def\TTT{{\mathbb T}}
\def\NNN{{\mathbb N}}

\def\de{{\rm d}}

\def\qed{\hfill\raise1pt\hbox{\vrule height5pt width5pt depth0pt}}

\def\ins#1#2#3{\vbox to0pt{\kern-#2 \hbox{\kern#1 #3}\vss}\nointerlineskip}

\begin{document}

\title{\bf Melnikov theory to all orders and Puiseux series\\
for subharmonic solutions}

\author{
\bf Livia Corsi and Guido Gentile
\vspace{2mm}
\\ \small 
Dipartimento di Matematica, Universit\`a di Roma Tre, Roma,
I-00146, Italy.
\\ \small
E-mail: lcorsi@mat.uniroma3.it, gentile@mat.uniroma3.it
}

\date{}

\maketitle

\begin{abstract}
We study the problem of subharmonic bifurcations for analytic systems
in the plane with perturbations depending periodically on time, in the
case in which we only assume that the subharmonic Melnikov function
has at least one zero. If the order of zero is odd, then there is
always at least one subharmonic solution, whereas if the order
is even in general other conditions have to be assumed to
guarantee the existence of subharmonic solutions. Even when such
solutions exist, in general they are not analytic in the perturbation
parameter. We show that they are analytic in a fractional power of the
perturbation parameter. To obtain a fully constructive algorithm
which allows us not only to prove existence but also to obtain
bounds on the radius of analyticity and to approximate the solutions
within any fixed accuracy, we need further assumptions.
The method we use to construct the solution -- when this is possible --
is based on a combination of the Newton-Puiseux algorithm
and the tree formalism. This leads to a graphical representation
of the solution in terms of diagrams. Finally, if the subharmonic
Melnikov function is identically zero, we show that it
is possible to introduce higher order generalisations,
for which the same kind of analysis can be carried out.
\end{abstract}

%%%%%%%%%%%%%%%%%%%%%%%%%%%%%%%%%%%%%%%%%%%%%%%%%%%%%%%%%%%%%%%%%%%%%%%%%
%%%%%%%%%%%%%%%%%%%%%%%%%%%%%%%%%%%%%%%%%%%%%%%%%%%%%%%%%%%%%%%%%%%%%%%%%
\zerarcounters
\section{Introduction}\label{sec:1}
%%%%%%%%%%%%%%%%%%%%%%%%%%%%%%%%%%%%%%%%%%%%%%%%%%%%%%%%%%%%%%%%%%%%%%%%%
%%%%%%%%%%%%%%%%%%%%%%%%%%%%%%%%%%%%%%%%%%%%%%%%%%%%%%%%%%%%%%%%%%%%%%%%%

\noindent The problem of subharmonic bifurcations was first considered by
Melnikov \cite{M}, who showed that the existence of
subharmonic solutions is related to the zeroes of a suitable function,
nowadays called the subharmonic Melnikov function.
The standard Melnikov theory usually studies the case
in which the subharmonic Melnikov function 
has a simple (i.e. first order) zero \cite{CH,GH}. In such a case
the problem can be reduced to a problem of implicit function theorem.

Nonetheless, it can happen that the subharmonic Melnikov function
either vanishes identically or has a zero which is of order higher
than one. In the first case hopefully one can go to higher orders,
and if a suitable higher order generalisation of
the subharmonic Melnikov function has a first order zero,
then one can proceed very closely to the standard case,
and existence of analytic subharmonic solutions is obtained.
Most of the papers in the literature consider this kind
of generalisations of Melnikov's theory, and often a second order
analysis is enough to settle the problem.

The second case is more subtle. The problem can be still
reduced to an implicit function problem, but the fact
that the zeroes are no longer simple prevents us from
applying the implicit function theorem. Thus, other arguments
must be used, based on the Weierstrass preparation theorem
and on the theory of the Puiseux series \cite{P,BK,CH,ALGM}.
However, a systematic analysis is missing in the literature.
Furthermore, in general, these arguments are not constructive:
if on the one hand they allow to prove (in certain cases)
the existence of at least one subharmonic solution,
on the other hand the problem of how many such solutions
really exist and how they can be explicitly constructed
has not been discussed in full generality.

The main difficulty for a constructive approach is that
the solution of the implicit function equation has to be
looked for by successive approximations. At each iteration step,
in order to find the correction to the approximate solution
found at the previous one, one has to solve a new implicit
function equation, which, in principle, still admits multiple roots.
So, as far as the roots of the equations are not simple,
one cannot give an algorithm to produce systematically
the corrections at the subsequent steps.

A careful discussion of a problem of the same kind can be
found in \cite{ALGM}, where the problem of bifurcations
from multiple limit cycles is considered -- cf. also \cite{P1,P2},
where the problem is further investigated. There, under the
hypothesis that a simple (real) zero is obtained at the first
iteration step, it is proved that the bifurcating solutions can be
expanded as fractional series (Puiseux series) of the
perturbation parameter. The method to compute the coefficients
of the series is based on the use of Newton's polygon \cite{BK,CH,ALGM},
and allows one to go to arbitrarily high orders. However,
the convergence of the series, and hence of the algorithm,
relies on abstract arguments of algebraic and geometric theory.

To the best of our knowledge, the case of subharmonic bifurcations 
was not discussed in the literature. Of course, in principle
one can think to adapt the same strategy as in \cite{ALGM}
for the bifurcations of limit cycles. But still, there are
issues which have not been discussed there. Moreover we have
a twofold aim. We are interested in results which are
both general -- not generic -- and constructive.
This means that we are interested in problems such as
the following one: which are the weaker conditions
to impose on the perturbation,
for a given integrable system and a given periodic solution,
in order to prove the existence of subharmonic solutions?
Of course the ideal result would be to have no restriction at all.
At the same time, we are also interested in explicitly construct
such solutions, within any prefixed accuracy.

The problem of subharmonic solutions in the case of multiple zeroes
of the Melnikov functions has been considered in \cite{ZL}, where
the following theorem is stated (without giving the proof) for $C^{r}$
smooth systems: if the subharmonic Melnikov function has a zero of
order $n\le r$, then there is at least one subharmonic solution. 
In any case the analyticity properties of the solutions are
not discussed. In particular the subharmonic solution is found as
a function of two parameters -- the perturbation parameter and the
initial phase of the solution to be continued --, but the
relation between the two parameters is not discussed.
We note that, in the analytic setting, it is exactly
this relation which produces the lack of analyticity
in the perturbation parameter. Furthermore, in \cite{ZL}
the case of zeroes of even order is not considered:
as we shall see, in that case the existence of
subharmonic solutions can not be proved in general,
but it can be obtained under extra assumptions.

In the remaining part of this section, we give a more detailed
account of our results.
One can formulate the problem both in the $C^{r}$ Whitney topology
and in the real-analytic setting. We shall choose the latter.
 From a technical point of view, this is mandatory since
our techniques requires for the systems to be analytic.
However, it is also very natural from a physical point of view,
because in practice in any physical applications the functions
appearing in the equations are analytic (often even polynomials),
and when they are not analytic they are not even smooth.
Also, we note since now that, even though we restrict our
analysis to the analytic setting, this does not mean at all
that we can not deal with problems where non-analytic phenomena
arise. The very case discussed in this paper provides a counterexample.

We shall consider systems which can be viewed as perturbations
of integrable systems, with the perturbation which depends
periodically in time. We shall use coordinates $(\a,A)$
such that, in the absence of the perturbation, $A$ is fixed
to a constant value, while $\a$ rotates on the circle:
hence all motions are periodic. As usual \cite{GH} we assume
that, for $A$ varying in a finite interval, the periods change
monotonically. Then we can write the equations of motion as
$\dot\a=\o(A)+\e F(\a,A,t)$, $\dot A = \e G(\a,A,t)$, with $G,F$
periodic in $\a$ and $t$. All functions are assumed to be analytic.
More formal definitions will be given in Section \ref{sec:2}.

Given a unperturbed periodic orbit $t\to (\a_{0}(t),A_{0}(t))$,
we define the subharmonic Melnikov function $M(t_{0})$ as the average
over a period of the function $G(\a_{0}(t),A_{0},t+t_{0})$.
By construction $M(t_{0})$ is periodic in $t_{0}$.
With the terminology introduced above, $\e$ is the perturbation
parameter and $t_{0}$ is the initial phase.
The following scenario arises.

\begin{itemize}

\item If $M(t_{0})$ has no zero, then there is no subharmonic solution,
that is no periodic solution which continues the unperturbed one
at $\e\neq0$.

\item Otherwise, if $M(t_{0})$ has zeroes, the following two cases
are possible: either $M(t_{0})$ has a zero of finite order
$\nnn$ or $M(t_{0})$ vanishes with
all its derivatives. In the second case, because of analyticity,
the function $M(t_{0})$ is identically zero.

\item If $M(t_{0})$ has a simple zero (i.e. $\nnn=1$),
then the usual Melnikov's theory applies. In particular
there exists at least one subharmonic
solution, and it is analytic in the perturbation parameter $\e$.

\item If $M(t_{0})$ has a zero of order $\nnn$, then in general
no result can be given about the existence of subharmonic
solutions. However one can introduce an infinite sequence
of polynomial equations, which are defined iteratively:
if the first equation admits a real non-zero root
and all the following equations admit a real root,
then a subharmonic solution exists, and it is
a function analytic in suitable fractional power of $\e$;
more precisely it is analytic in $\h=\e^{1/p}$,
for some $p\le \nnn!$, and hence it is analytic in $\e^{1/\nnn!}$.
If at some step the root is simple,
an algorithm can be given in order to construct
recursively all the coefficients of the series.

\item If we further assume that the order $\nnn$ of the zero is odd,
then we have that all the equations of the sequence satisfy
the request made above on the roots, so that we can conclude
that in such a case at least one subharmonic solution exists.
Again, in order to really construct the solution, by providing
an explicit recursive algorithm, we need that at a certain level
of the iteration scheme a simple root appears.

\item Moreover we have at most $\nnn$ periodic solutions bifurcating
from the unperturbed one with initial phase $t_{0}$.
Of course, to count all subharmonic solutions we have also
to sum over all the zeroes of the subharmonic Melnikov function.

\item Finally, if $M(t_{0})$ vanishes identically as a function of
$t_{0}$, then we have to extend the analysis up to second order,
and all the cases discussed above for $M(t_{0})$ have to
repeated for a suitable function $M_{1}(t_{0})$, which is obtained
in the following way. If $M(t_{0})\equiv 0$ then the solution
$t \to (\a(t),A(t))$ is defined up to first order -- as it is
easy to check --, so that one can expand the function
$G(\a(t),A(t),t+t_{0})$ up to first order: we call $M_{1}(t_{0})$
its average over a period of the unperturbed solution.
In particular if also $M_{1}(t_{0})$ vanishes identically then
one can push the perturbation theory up to second order,
and, after expanding the function $G(\a(t),A(t),t+t_{0})$
up to second order, one defines $M_{2}(t_{0})$ as its average
over a period, and so on.

\end{itemize}

The first conclusion we can draw is that in general we cannot
say that for any vector field $(F,G)$ there is at least
one subharmonic solution of given period. We need some condition on $G$.
We can require for $G$ to be a zero-mean function, so that
it has at least one zero of odd order. For instance, this holds true
if the vector field is Hamiltonian, since in such a case $G$
is the $\alpha$-derivative of a suitable function. The same result
follows if the equations describe a Hamiltonian system in the presence
of small friction -- how small depends on the particular
resonance one is looking at \cite{GBD}. But of course,
all these conditions are stronger than what is really needed.

A second conclusion is that, even when a subharmonic solution
turns out to exist (and to be analytic in a suitable
fractionary power of the perturbation parameter), a constructive
algortithm to compute it within any given accuracy cannot be
provided in general. This becomes possible only if some further
assumption is made. So there are situations where
one can obtain an existence result of the solution,
but the solutin itself cannot be constructed.
Note that such situations are highly non-generic, because
they arise if one finds at each iterative step a polynomial
with multiple roots -- which is a non-generic case;
cf. Appendix \ref{app:A}.

The methods we shall use to prove the results above will be
of two different types. We shall rely on standard general techniques,
based on the Weierstrass preparation theorem, in order to show
that under suitable assumptions the solutions exist and to prove
in this case the convergence of the series.
Moreover, we shall use a combination of the Newton-Puiseux process
and the diagrammatic techniques based on the tree formalism
\cite{GGG,GBD1,GBD} in order to provide a recursive algorithm,
when possible. Note that in such a case the convergence
of the Puiseux series follows by explicit construction of the
coefficients, and an explicit bound of the
radius of convergence is obtained through the estimates of
the coefficients -- on the contrary there is no way to provide
quantitative bounds with the aforementioned abstract arguments.
These results extend those in \cite{GBD}, where a special case
was considered.

The paper is organised as follows. In Section \ref{sec:2} we
formulate rigorously the problem of subharmonic bifurcations
for analytic ordinary differential equations in the plane,
and show that, if the subharmonic Melnikov function admits
a finite order zero, the problem can be reduced to an analytic
implicit equation problem -- analyticity will be proved
in Section \ref{sec:5} by using the tree formalism.
In Section \ref{sec:3} we discuss
the Newton-Puiseux process, which will be used to iteratively
attack the problem. At each iteration step one has to solve
a polynomial equation. Thus, in the complex setting \cite{BK}
the process can be pushed forward indefinitely, whereas
in the real setting one has to impose at each step that a
real root exists. If the order of zero of the
subharmonic Melnikov function is odd, the latter condition
is automatically satisfied, and hence the existence of at least
one subharmonic solution is obtained (Theorem \ref{thm:1}).
If at some step of the iteration a simple root appears,
then we can give a fully constructive algorithm which allows us
to estimate the radius of analyticity and to approximate
the solution within any fixed accuracy (Theorem \ref{thm:2}).
This second result will be proved in Sections \ref{sec:4}
and \ref{sec:6}, again by relying on the tree formalism;
some more technical aspects of the proof will be dealt with
in Appendix \ref{app:B}. Finally in Section \ref{sec:7}
we consider the case in which the subharmonic Melnikov function
vanishes identically, so that one has to repeat the analysis
for suitable higher order generalisations of that function.
This will lead to Theorems \ref{thm:3} and \ref{thm:4},
which generalise Theorems \ref{thm:2} and \ref{thm:1}, respectively.

%%%%%%%%%%%%%%%%%%%%%%%%%%%%%%%%%%%%%%%%%%%%%%%%%%%%%%%%%%%%%%%%%%%%%%%%%
%%%%%%%%%%%%%%%%%%%%%%%%%%%%%%%%%%%%%%%%%%%%%%%%%%%%%%%%%%%%%%%%%%%%%%%%%
\zerarcounters
\section{Set-up}\label{sec:2}
%%%%%%%%%%%%%%%%%%%%%%%%%%%%%%%%%%%%%%%%%%%%%%%%%%%%%%%%%%%%%%%%%%%%%%%%%
%%%%%%%%%%%%%%%%%%%%%%%%%%%%%%%%%%%%%%%%%%%%%%%%%%%%%%%%%%%%%%%%%%%%%%%%%

\noindent Let us consider the ordinary differential equation
\be\left\{\ba
&\dot{\a}=\o(A)+\e F(\a,A,t), \\
&\dot{A}=\e G(\a,A,t),
\label{eq:2.1}\ea\right.\ee
where $(\a,A)\in\MM:=\TTT\times W$, with $W\subset\RRR$ an open set,
the map $A\mapsto\o(A)$ is real analytic in $A$, and the functions
$F,\,G$ depend analytically on their arguments and are $2\p$-periodic
in $\a$ and $t$. Finally $\e$ is a real parameter.

Set $\a_{0}(t)=\o(A_{0})t$ and $A_{0}(t)=A_{0}$. In the
\emph{extended phase-space} $\MM\times\RRR$, for $\e=0$, the solution
$(\a_{0}(t),A_{0}(t),t+t_{0})$ describes an invariant torus, which is
uniquely determined by the ``energy'' $A_{0}$. Hence the motion of the
variables $(\a,A,t)$ is quasi-periodic, and reduces to a periodic
motion whenever $\o(A_{0})$ becomes commensurate with $1$. If $\o(A_{0})$ is
rational we say that the torus is \emph{resonant}. The parameter
$t_{0}$ will be called the \emph{initial phase}: it fixes the initial
datum on the torus. Only for some values of the parameter $t_{0}$
periodic solutions lying on the torus are expected to persist
under perturbation: such solutions are called \emph{subharmonic solutions}.

Denote by $T_{0}(A_{0})=2\p/\o(A_{0})$ the period of the trajectories on
the unperturbed torus, and define $\o'(A):={\de}\o(A)/{\de}A$.
If $\o(A_{0})=p/q\in\QQQ$, call $T=T(A_{0})=2\p q$ the period
of the trajectories in the extended phase space. We shall call $p/q$
the \emph{order} of the corresponding subharmonic solutions.

%%%%%%%%%%%%%%%%%%%%%%%%%%%%%%%%%%%%%%%%%%%%%%%%%%%%%%%%%%%%%%%%%%%%%%%%%
\begin{hyp}\label{hyp1} One has $\o'(A_{0})\neq0$.\end{hyp}
%%%%%%%%%%%%%%%%%%%%%%%%%%%%%%%%%%%%%%%%%%%%%%%%%%%%%%%%%%%%%%%%%%%%%%%%%

Define
\be
M(t_{0}) := \frac{1}{T}\int_{0}^{T} \de t\, G(\a_{0}(t),A_{0},t+t_{0}),
\label{eq:2.2}\ee
which is called the \emph{subharmonic Melnikov function of order} $q/p$.
Note that $M(t_{0})$ is $2\p$-periodic in $t_{0}$.

%%%%%%%%%%%%%%%%%%%%%%%%%%%%%%%%%%%%%%%%%%%%%%%%%%%%%%%%%%%%%%%%%%%%%%%%%
\begin{hyp}\label{hyp2} There exist $t_{0}\in[0,2\p)$ and $\nnn\in\NNN$
such that
\be
\frac{\de^{k}}{\de t_{0}^{k}}M(t_{0})=0\;\;\;\forall\;0
\leq k\leq\nnn-1,\phantom{ and }
D(t_{0}) := \frac{\de^{\nnn}}{\de t_{0}^{\nnn}}M(t_{0})\neq0 ,
\label{eq:2.3}\ee
that is $t_{0}$ is a zero of order $\nnn$ for
the subharmonic Melnikov function.
\end{hyp}
%%%%%%%%%%%%%%%%%%%%%%%%%%%%%%%%%%%%%%%%%%%%%%%%%%%%%%%%%%%%%%%%%%%%%%%%%

For notational semplicity, we shall not make explicit the dependence
on $t_{0}$ most of times; for instance we shall write $D(t_{0})=D$.
For any $T$-periodic function $F$ we shall denote by $\mean{F}$
its average over the period $T$.

The solution of (\ref{eq:2.1}) with initial conditions
$(\a(0),A(0))$ can be written as
\begin{equation}
\left( \begin{matrix} \a(t) \\
A(t) \end{matrix} \right)
= W(t) \left( \begin{matrix} \a(0) \\
A(0) \end{matrix} \right) + W(t) \int_{0}^{t} {\rm d} \tau \,
W^{-1}(\tau) \left( \begin{matrix} \Phi(\tau) \\
\Gamma(\tau) \end{matrix} \right) ,
\label{eq:2.4} \end{equation}
where we have denoted by
\begin{equation}
W(t) = \left( \begin{matrix}
1 & \o'(A_{0}) t \\ 0 & 1 \end{matrix} \right)
\label{eq:2.5} \end{equation}
the Wronskian matrix solving the linearised system, and set
\begin{equation}
\Phi(t) = \e F(t) + \o(A(t)) - \o(A_{0}) -
\o'(A_{0}) \left( A(t) - A_{0} \right) ,
\qquad \Gamma(t) = \e G(t) .
\label{eq:2.6} \end{equation}
shortening $F(t)=F(\a(t),A(t),t+t_{0})$
and $G(t)=G(\a(t),A(t),t+t_{0})$.

By using explicitly (\ref{eq:2.5}) in (\ref{eq:2.4}) we obtain
\begin{equation}
\begin{cases}
{\displaystyle
\a(t) = \a(0) + t \, \o'(A_{0})\,A(0) +
\int_{0}^{t} {\rm d}\tau \, \Phi(\tau) + \o'(A_{0})
\int^{t}_{0} {\rm d} \tau \int_{0}^{\tau} {\rm d}\tau'
\Gamma(\tau') , } \\
{\displaystyle A(t) = A(0) +
\int_{0}^{t} {\rm d}\tau \, \Gamma(\tau) , }
\end{cases}
\label{eq:2.7} \end{equation}
with the notations (\ref{eq:2.6}).

In order to obtain a periodic solution we need for the
mean $\langle \Gamma \rangle$ of the function $\Gamma$ to be zero.
In this case, if we fix also
\begin{equation}
A(0) = - \frac{1}{\ \o'(A_{0})} \langle \Phi \rangle -
\langle \GG \rangle , \qquad
\GG(\tau) = \int_{0}^{\tau} {\rm d}\tau' \left(
\Gamma(\tau') - \langle \Gamma \rangle \right) ,
\label{eq:2.8} \end{equation}
then the corresponding solution turns out to be periodic.
So, instead of (\ref{eq:2.7}), we consider the system
\begin{equation}
\begin{cases}
{\displaystyle
\a(t) = \a(0) + \int_{0}^{t} {\rm d}\tau \,
\left( \Phi(\tau) -\langle \Phi \rangle \right) + \o'(A_{0})
\int^{t}_{0} {\rm d} \tau (\GG(\t)-\mean{\GG}), } \\
{\displaystyle A(t) = A(0) + \GG(t) , } \\
{\displaystyle \langle \Gamma \rangle = 0 , }
\end{cases}
\label{eq:2.9} \end{equation}
where $A(0)$ is determined according to (\ref{eq:2.8})
and $\a(0)$ is considered as a free parameter.

We start by considering the \textit{auxiliary system}
\begin{equation}
\begin{cases}
{\displaystyle
\a(t) = \a(0) + \int_{0}^{t} {\rm d}\tau \,
\left( \Phi(\tau) -\langle \Phi \rangle \right) + \o'(A_{0})
\int^{t}_{0} {\rm d} \tau (\GG(\t)-\mean{\GG}), } \\
{\displaystyle A(t) = A(0) + \GG(t) , } \\
\end{cases}
\label{eq:2.10} \end{equation}
that is we neglect for the moment the condition that the mean
of $\Gamma$ has to be zero. Of course, only in that case the
solution of (\ref{eq:2.10}) is solution also of (\ref{eq:2.9}),
hence of (\ref{eq:2.7}).

It can be more convenient to work in Fourier space.
As we are looking for periodic solutions of period $T=2\pi q$,
i.e. of frequency $\o=1/q$, we can write
\begin{equation}
\a(t) = \a_{0}(t) + \beta(t) , \qquad
\beta(t) = \sum_{\n\in\ZZZ} {\rm e}^{i\o\n t} \beta_{\n} , \qquad
A(t) = A_{0} + B(t) , \qquad
B(t) = \sum_{\n\in\ZZZ} {\rm e}^{i\o\n t} B_{\n} .
\label{eq:2.11} \end{equation}
If we expand
\begin{equation}
G(\a,A,t+t_{0}) = \sum_{\nu\in\ZZZ} \sum_{\nu'\in\ZZZ}
{\rm e}^{i\nu\a + i\nu' (t+t_{0})}
G_{\nu,\nu'}(A) , \qquad G_{\nu,\nu'}(A,t_{0}) :=
{\rm e}^{i\nu' t_{0}} G_{\nu,\nu'}(A) ,
\label{eq:2.12} \end{equation}
with an analogous expressions for the function $\Phi(t)$,
then we can write
\begin{equation}
\Gamma(t) = \sum_{\n\in\ZZZ} {\rm e}^{i\o\n t} \Gamma_{\n} , \qquad
\Phi(t) = \sum_{\n\in\ZZZ} {\rm e}^{i\o\n t} \Phi_{\n} ,
\label{eq:2.13} \end{equation}
with
\begin{subequations}
\begin{align}
\Gamma_{\nu} & =
\e \sum_{r=0}^{\io} \sum_{s=0}^{\io}
\sum_{p\nu_{0}+q\nu_{0}' + \nu_{1} + \ldots + \nu_{r+s}= \nu}
\frac{1}{r!s!} (i\nu_{0})^{r} \partial_{A}^{s}
G_{\nu_{0},\nu_{0}'}(A_{0},t_{0}) \,
\beta_{\nu_{1}} \ldots \beta_{\nu_{r}}
B_{\nu_{r+1}} \ldots B_{\nu_{r+s}} ,
\label{eq:2.14a} \\
\Phi_{\nu} & =
\e \sum_{r=0}^{\io} \sum_{s=0}^{\io}
\sum_{p\nu_{0}+q\nu_{0}' + \nu_{1} + \ldots + \nu_{r+s}= \nu}
\frac{1}{r!s!} (i\nu_{0})^{r} \partial_{A}^{s}
F_{\nu_{0},\nu_{0}'}(A_{0},t_{0}) \,
\beta_{\nu_{1}} \ldots \beta_{\nu_{r}}
B_{\nu_{r+1}} \ldots B_{\nu_{r+s}} \nonumber \\
& + \sum_{s=2}^{\io} \sum_{\nu_{1} + \ldots + \nu_{s}= \nu}
\frac{1}{s!} \partial_{A}^{s}
\o(A_{0})\,B_{\nu_{1}} \ldots B_{\nu_{s}} .
\label{eq:2.14b}
\end{align}
\label{eq:2.14}
\end{subequations}
\vskip-.3truecm
\noindent Then (\ref{eq:2.10}) becomes
\begin{equation}
\begin{cases}
{\displaystyle
\beta_{\n} = \frac{\Phi_{\n}}{i\o\n} +
\o'(A_{0}) \frac{\Gamma_{\n}}{(i\o\n)^{2}} } , \\
{\displaystyle B_{\n} = \frac{\Gamma_{\n}}{i\o\n} } ,
\end{cases}
\label{eq:2.15} \end{equation}
for $\n\neq0$, provided
\begin{equation}
\begin{cases}
{\displaystyle
\beta_{0} = \a(0) -
\sum_{\substack{\n\in\ZZZ \\ \n\neq 0}} \frac{\Phi_{\n}}{i\o\n} -
\o'(A_{0}) \sum_{\substack{\n\in\ZZZ \\ \n\neq 0}}
\frac{\Gamma_{\n}}{(i\o\n)^{2}} } , \\
{\displaystyle B_{0} = A(0) - \sum_{\substack{\n\in\ZZZ \\ \n\neq 0}}
\frac{\Gamma_{\n}}{i\o\n} = - \frac{\Phi_{0}}{\o'(A_{0})} } ,
\end{cases}
\label{eq:2.16} \end{equation}
for $\n=0$. Also (\ref{eq:2.9}) can be written in the same form, with
the further constraint $\Gamma_{0}=0$.

Then we can use $\beta_{0}$ as a free parameter, instead of $\a(0)$.
This means that we look for a value of $\beta_{0}$
(depending on $\e)$ such that, by defining $B_{0}$ according
to the second equation in (\ref{eq:2.16}), the coefficients
$\beta_{\n},B_{\n}$ are given by (\ref{eq:2.15}) for $\n\neq0$.
In other words, in Fourier space (\ref{eq:2.10}) becomes
\begin{equation}
\begin{cases}
{\displaystyle
\beta_{\n} = \frac{\Phi_{\n}}{i\o\n} +
\o'(A_{0}) \frac{\Gamma_{\n}}{(i\o\n)^{2}} } , \qquad
{\displaystyle B_{\n} = \frac{\Gamma_{\n}}{i\o\n} } ,
\qquad \n\neq0 , \\
{\displaystyle B_{0} = - \frac{\Phi_{0}}{\o'(A_{0})} } ,
\end{cases}
\label{eq:2.17} \end{equation}
whereas $\beta_{0}$ is left as a free parameter.

We look for a solution $(\ol{\a}(t),\ol{A}(t))$ of (\ref{eq:2.10})
which can be written
as a formal Taylor series in $\e$ and $\beta_{0}$, so that
\begin{subequations}
\begin{align}
\ol{\a}(t) & = 
\ol{\a}(t;\e,\beta_{0}) = \a_{0}(t) +
\sum_{k=1}^{\infty} \sum_{j=0}^{\infty} \e^{k} \beta_{0}^{j} \, 
\ol{\b}^{(k,j)}(t) ,
\label{eq:2.18a} \\
\ol{A}(t) & =
\ol{A}(t;\e,\beta_{0}) = A_{0} + \sum_{k=1}^{\infty}
\sum_{j=0}^{\infty} \e^{k} \beta_{0}^{j} \,
\ol{B}^{(k,j)}(t) ,
\label{eq:2.18b}
\end{align}
\label{eq:2.18} \end{subequations}
\vskip-.3truecm
\noindent which reduces to $(\a_{0}(t),A_{0})$ as $\e\to0$.
By comparing (\ref{eq:2.18}) with (\ref{eq:2.11}) we can write
the Fourier coefficients of the solution $(\ol{\a}(t),\ol{A}(t))$,
for $\n\neq0$, as
\begin{equation}
\ol{\b}_{\n} = \ol{\b}_{\n}(\e,\beta_{0}) = \sum_{k=1}^{\infty}
\sum_{j=0}^{\infty} \e^{k} \beta_{0}^{j} \, \ol{\b}^{(k,j)}_{\n} ,
\qquad \ol{B}_{\n} = \ol{B}_{\n}(\e,\beta_{0}) = \sum_{k=1}^{\infty}
\sum_{j=0}^{\infty} \e^{k} \beta_{0}^{j} \, \ol{B}^{(k,j)}_{\nu} ,
\label{eq:2.19} \end{equation}
and, analogously, $\ol{B}_{0}^{(k,j)}$ is the contribution
to order $k$ in $\e$ and $j$ in $\beta_{0}$ to $\ol{B}_{0}$.

By analyticity also the function $\ol{\Gamma}(t)=
\e G(\ol{\a}(t),\ol{A}(t),t+t_{0})$ can be formally
expanded in powers of $\e$ and $\beta_{0}$, and one has
\begin{equation}
\ol{\Gamma}(t) = \ol{\Gamma}(t;\e,\beta_{0}) =
\sum_{k=1}^{\infty} \sum_{j=0}^{\io} \e^{k} \beta_{0}^{j}
\ol{\Gamma}^{(k,j)}(t) = \sum_{k=1}^{\infty}
\sum_{j=0}^{\io} \e^{k} \beta_{0}^{j}
\sum_{\n\in\ZZZ} {\rm e}^{i\o\n t} \ol{\Gamma}^{(k,j)}_{\n} ,
\label{eq:2.20} \end{equation}
where each $\ol{\Gamma}^{(k,j)}_{\n}$ is expressed in terms of
the Taylor coefficients of (\ref{eq:2.19}) of order strictly less
than $k,j$. By definition one has $\langle \ol{\Gamma}^{(k,j)}
\rangle = \ol{\Gamma}^{(k,j)}_{0}$. The same considerations hold for
$\ol{\Phi}(t)=\Phi(\ol{\a}(t),\ol{A}(t),t+t_{0})$.

Hence one can formally write, for all $k\geq1$ and $j\geq0$
\be
\left\{\ba
&\ol{\b}^{(k,j)}_{\n} = \frac{\ol{\F}^{(k,j)}_{\n}}{i\o\n}+
\o'(A_{0}) \frac{\ol{\G}^{(k,j)}_{\n}}{(i\o\n)^2}, \qquad
\ol{B}^{(k,j)}_{\n} = \frac{\ol{\G}^{(k,j)}_{\n}}{i\o\n},
\qquad \n\neq0 \\
&\ol{B}^{(k,j)}_{0}=- \frac{\ol{\F}^{(k,j)}_{0}}{\o'(A_{0})}.
\ea\right.
\label{eq:2.21} \ee

%%%%%%%%%%%%%%%%%%%%%%%%%%%%%%%%%%%%%%%%%%%%%%%%%%%%%%%%%%%%%%%%%%%%%%%%%
\begin{lemma} \label{lem:1}
For any $\beta_{0}\in\RRR$ the system (\ref{eq:2.10}) admits a
solution $(\ol{\a}(t),\ol{A}(t))$ which is $T$-periodic
in time and analytic in $\e$, depending analytically
on the parameter $\beta_{0}$.
\end{lemma}
%%%%%%%%%%%%%%%%%%%%%%%%%%%%%%%%%%%%%%%%%%%%%%%%%%%%%%%%%%%%%%%%%%%%%%%%%

%%%%%%%%%%%%%%%%%%%%%%%%%%%%%%%%%%%%%%%%%%%%%%%%%%%%%%%%%%%%%%%%%%%%%%%%%
\prova One can use the tree formalism introduced in Section \ref{sec:5};
see in particular Proposition \ref{prop:1}.\qed
%%%%%%%%%%%%%%%%%%%%%%%%%%%%%%%%%%%%%%%%%%%%%%%%%%%%%%%%%%%%%%%%%%%%%%%%%
\salto
%%%%%%%%%%%%%%%%%%%%%%%%%%%%%%%%%%%%%%%%%%%%%%%%%%%%%%%%%%%%%%%%%%%%%%%%%

It can be convenient to introduce also the Taylor coefficients
\be\label{eq:2.22}\ba
&{\eknu{\ol{\b}} k \n}(\b_{0})=
\sum_{j\geq0}\b_{0}^{j}\;{\eknu{\ol{\b}} {k,j} \n},\qquad
{\eknu{\ol{B}} k \n}(\b_{0})=
\sum_{j\geq0}\b_{0}^{j}\;{\eknu{\ol{B}} {k,j} \n},\\
&{\eknu{\ol{\G}} k \n}(\b_{0})=
\sum_{j\geq0}\b_{0}^{j}\;{\eknu{\ol{\G}} {k,j} \n},\qquad
{\eknu{\ol{\F}} k \n}(\b_{0})=
\sum_{j\geq0}\b_{0}^{j}\;{\eknu{\ol{\F}} {k,j} \n},
\ea\ee
Note that $\ol{\Gamma}^{(k)}_{\nu}(0)=\ol{\Gamma}^{(k,0)}_{\nu}$,
and so on.

%%%%%%%%%%%%%%%%%%%%%%%%%%%%%%%%%%%%%%%%%%%%%%%%%%%%%%%%%%%%%%%%%%%%%%%%%
\begin{lemma} \label{lem:2}
Consider the system (\ref{eq:2.10}).
Assume that $\ol{\Gamma}^{(k)}_{0}(0) = 0$ for all $k\in\NNN$.
Then for $\beta_{0}=0$ the solution $(\ol{\a}(t),\ol{A}(t))$
of (\ref{eq:2.10}) is also a solution of (\ref{eq:2.9}).
\end{lemma}
%%%%%%%%%%%%%%%%%%%%%%%%%%%%%%%%%%%%%%%%%%%%%%%%%%%%%%%%%%%%%%%%%%%%%%%%%

%%%%%%%%%%%%%%%%%%%%%%%%%%%%%%%%%%%%%%%%%%%%%%%%%%%%%%%%%%%%%%%%%%%%%%%%%
\prova Simply note that (\ref{eq:2.10}) reduces to (\ref{eq:2.9})
if $\ol{\Gamma}^{(k)}_{0}(0) = 0$ for all $k\in\NNN$. \qed
%%%%%%%%%%%%%%%%%%%%%%%%%%%%%%%%%%%%%%%%%%%%%%%%%%%%%%%%%%%%%%%%%%%%%%%%%
\salto
%%%%%%%%%%%%%%%%%%%%%%%%%%%%%%%%%%%%%%%%%%%%%%%%%%%%%%%%%%%%%%%%%%%%%%%%%

Of course we expect in general that $\ol{\Gamma}^{(k)}_{0}(0)$ do not
vanish for all $k\in\NNN$. In that case, let $k_{0}\in\NNN$ be such that
$\ol{\Gamma}^{(k)}_{0}(0)=0 $ for $k=1,\ldots,k_{0}$
and $\ol{\Gamma}^{(k_{0}+1)}_{0}(0) \neq 0$.

Let us define
\be\label{eq:2.23}
\FF^{(0)}(\e,\b_{0}):=\sum_{k,j\geq0}\e^{k}\b_{0}^{j}\FF_{k,j}^{(0)},
\qquad
\FF_{k,j}^{(0)}=\ol{\G}^{(k+1,j)}_{0},
\ee
so that $\e\FF^{(0)}(\e,\b_{0})=\mean{\ol{\G}(\;\cdot\;;\e,\b_{0})}$.

%%%%%%%%%%%%%%%%%%%%%%%%%%%%%%%%%%%%%%%%%%%%%%%%%%%%%%%%%%%%%%%%%%%%%%%%%
\begin{lemma} \label{lem:3}
$\FF^{(0)}(\e,\b_{0})$ is $\b_{0}$-general of order $\nnn$, \ie
$\FF_{0,j}^{(0)}=0$
for $j=0,\ldots,\nnn-1$, while $\FF_{0,\nnn}^{(0)}\neq0$.
\end{lemma}
%%%%%%%%%%%%%%%%%%%%%%%%%%%%%%%%%%%%%%%%%%%%%%%%%%%%%%%%%%%%%%%%%%%%%%%%%

%%%%%%%%%%%%%%%%%%%%%%%%%%%%%%%%%%%%%%%%%%%%%%%%%%%%%%%%%%%%%%%%%%%%%%%%%
\prova This can be easily shown using the tree formalism
introduced in Section \ref{sec:4}.
In fact for all $j$, $\FF_{0,j}^{(0)}=\ol{\G}^{(1,j)}_{0}$ is associated
with a tree with $1$ node and $j$ leaves. Hence one has
\be\label{eq:2.24}
j!\ol{\G}^{(1,j)}_{0}= \mean{\dpr_{\a}^{j}G(\a_{0}(\cdot),A_{0},
\cdot+t_{0})}
= (-\o(A_{0}))^{-j} \frac{\de^{j}M}{\de t_{0}^{j}}(t_{0}),
\ee
where the second equality is provided by Lemma 3.9 on \cite{GBD}.
Then $\FF^{(0)}(\e,\b_{0})$ is $\b_{0}$-general of order $\nnn$
by Hypothesis \ref{hyp2}.\qed
%%%%%%%%%%%%%%%%%%%%%%%%%%%%%%%%%%%%%%%%%%%%%%%%%%%%%%%%%%%%%%%%%%%%%%%%%
\salto
%%%%%%%%%%%%%%%%%%%%%%%%%%%%%%%%%%%%%%%%%%%%%%%%%%%%%%%%%%%%%%%%%%%%%%%%%

Our aim is to find $\b_{0}=\b_{0}(\e)$ such that $\FF^{(0)}(\e,\b_{0}(\e))
\equiv0$. For such $\b_{0}$ a solution of (\ref{eq:2.10}) is also solution
of (\ref{eq:2.9}). If we are successful in doing so,
then we have proved the existence of subharmonic solutions.

%%%%%%%%%%%%%%%%%%%%%%%%%%%%%%%%%%%%%%%%%%%%%%%%%%%%%%%%%%%%%%%%%%%%%%%%%
%%%%%%%%%%%%%%%%%%%%%%%%%%%%%%%%%%%%%%%%%%%%%%%%%%%%%%%%%%%%%%%%%%%%%%%%%
\zerarcounters
\section{The Newton-Puiseux process and main results}
\label{sec:3}
%%%%%%%%%%%%%%%%%%%%%%%%%%%%%%%%%%%%%%%%%%%%%%%%%%%%%%%%%%%%%%%%%%%%%%%%%
%%%%%%%%%%%%%%%%%%%%%%%%%%%%%%%%%%%%%%%%%%%%%%%%%%%%%%%%%%%%%%%%%%%%%%%%%

\noindent Given a convergent power series $\FF^{(0)}(\e,\b_{0})
\in\RRR\{\e,\b_{0}\}$ as in (\ref{eq:2.23}),
we call \emph{carrier of} $\FF^{(0)}$ the set
\be\label{eq:3.1}
\D(\FF^{(0)}) := \{(k,j)\in \NNN\times\NNN\,:\, \FF_{k,j}^{(0)}\neq0\}.
\ee
For all $v\in\D(\FF^{(0)})$ let us consider the positive quadrant
$\AAA_v:=\{v\}+(\RRR_+)^2$ moved up to $v$, and
define
\be\label{eq:3.2}\AAA:=\bigcup_{v\in\D(\FF^{(0)})}\AAA_v.\ee
Let $\CC$ be the convex hull of $\AAA$. The boundary $\dpr\CC$ consists
of a compact polygonal path
$\PP^{(0)}$ and two half lines $\RR_{1}^{(0)}$ and $\RR_{2}^{(0)}$.
The polygonal path $\PP^{(0)}$ is called the
\emph{Newton polygon} of $\FF^{(0)}$.

\begin{figure}[!ht]
\begin{center}
{
\psfrag{n}{$\nnn$}\psfrag{j}{$j$}
\psfrag{R1}{$\RR^{(0)}_{1}$}\psfrag{R2}{$\RR^{(0)}_{2}$}
\psfrag{k}{$k$}\psfrag{k0}{$k_{0}$}
\psfrag{Pi}{$\PP^{(0)}_{i}$}\psfrag{ri}{$\rrr^{(0)}_{i}$}
\psfrag{rimi}{$\frac{\rrr^{(0)}_{i}}{\m^{(0)}_{i}}$}
\includegraphics[width=10cm]{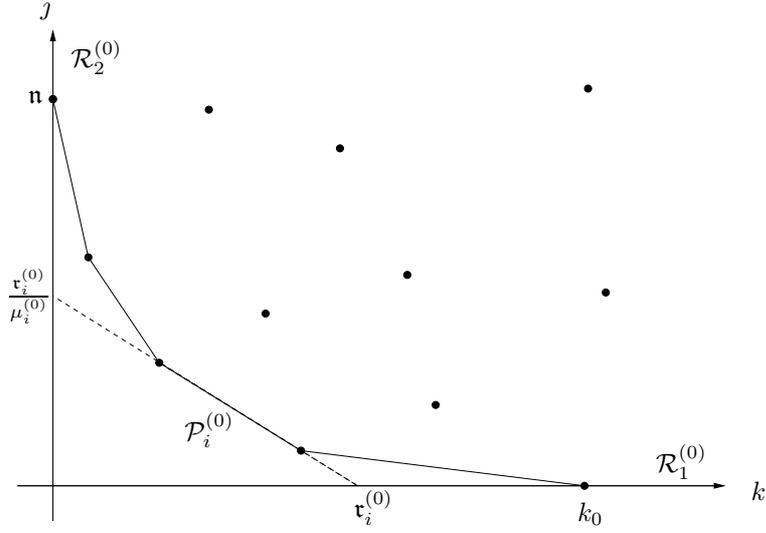}
}
\caption{\label{fig:1}
\footnotesize{Newton polygon.}}
\end{center}
\end{figure}

Notice that if the Newton polygon is a single point or,
more generally, if $\FF_{k,0}^{(0)}=0$ for all
$k\geq0$ then there exists $\ol{\j}\geq1$ such that
$\FF^{(0)}(\e,\b_{0})=\b_{0}^{\ol{\j}}\cdot\ol{\GG}(\e,\b_{0})$
with $\ol{\GG}(\e,0)\neq0$, hence $\b_{0}\equiv0$ is 
a solution of equation $\FF^{(0)}(\e,0)=0$, that is the
conclusion of Lemma \ref{lem:2}.
Otherwise, if we further assume that $\FF^{(0)}$ is $\b_0$-general
of some finite order $\nnn$,
there is at least a point of $\D(\FF^{(0)})$ on each axis,
then the Newton polygon $\PP^{(0)}$ is formed by $N_{0}\geq1$ segments
$\PP_{1}^{(0)},\ldots,\PP_{N_{0}}^{(0)}$ and we write $\PP^{(0)}=
\PP_{1}^{(0)}\cup\ldots\cup\PP_{N_{0}}^{(0)}$; cf. Figure \ref{fig:1}.

%%%%%%%%%%%%%%%%%%%%%%%%%%%%%%%%%%%%%%%%%%%%%%%%%%%%%%%%%%%%%%%%%%%%%%%%%
% figure 1
%%%%%%%%%%%%%%%%%%%%%%%%%%%%%%%%%%%%%%%%%%%%%%%%%%%%%%%%%%%%%%%%%%%%%%%%%
%%%%%%%%%%%%%%%%%%%%%%%%%%%%%%%%%%%%%%%%%%%%%%%%%%%%%%%%%%%%%%%%%%%%%%%%%

For all $i=1,\ldots,N_{0}$ let $-1/\m_{i}^{(0)}\in\QQQ$ be the slope
of the segment $\PP_{i}^{(0)}$, so that one can partition $\FF^{(0)}$
according to the weights given by $\m_{i}^{(0)}$:
\be
\FF^{(0)}(\e,\b_{0})=\tilde{\FF}_{i}^{(0)}(\e,\b_{0}) +
\GG_{i}^{(0)}(\e,\b_{0}) =
\!\!\!\!\!\!
\sum_{k+j\m_{i}^{(0)} = \rrr_{i}^{(0)}}
\!\!\!\!\!\!
\FF_{k,j}^{(0)}\e^{k}\b_{0}^{j} \,\, +
\!\!\!\!\!\!
\sum_{k+j\m_{i}^{(0)} > \rrr_{i}^{(0)}}
\!\!\!\!\!\!
\FF_{k,j}^{(0)}\e^{k}\b_{0}^{j},
\label{eq:3.3}\ee
where $\rrr_{i}^{(0)}$ is the intercept on the $k$-axis of the
continuation of $\PP_{i}^{(0)}$.

Hence the first approximate solutions of $\FF(\e,\b_{0})=0$
are the solutions of the quasi-homogeneous equations
\be\label{eq:3.4}
\tilde{\FF}_{i}^{(0)} (\e,\b_{0}) =
\!\!\!\!\!\!\!\!
\sum_{k\ppp_{i}^{(0)}+j\hhh_{i}^{(0)}= \sss_{i}^{(0)}}
\!\!\!\!\!\!\!\!
\FF_{k,j}^{(0)}\e^{k}\b_{0}^{j}=0, \qquad i=1,\ldots,N_{0},
\ee
where $\hhh_{i}^{(0)}/\ppp_{i}^{(0)}=\m_{i}^{(0)}$, with
$\hhh_{i}^{(0)}, \ppp_{i}^{(0)}$ relatively prime integers,
and $\sss_{i}^{(0)} = \ppp_{i}^{(0)} \rrr_{i}^{(0)}$.

We introduce the polynomials ${\eknu P 0 i} = {\eknu P 0 i}(c)$
in such a way that
\be\label{eq:3.5}
\tilde{\FF}_{i}^{(0)} (\e,c(\s_{0}\e)^{\m_{i}^{(0)}}) =
(\s_{0}\e)^{\rrr_{i}^{(0)}}
\!\!\!\!\!\!\!\!
\sum_{k\ppp_{i}^{(0)}+j\hhh_{i}^{(0)}= \sss_{i}^{(0)}}
\!\!\!\!\!\!\!\!
Q_{k,j}c^{j}=(\s_{0}\e)^{\rrr_{i}^{(0)}}
{\eknu P 0 i}(c), \qquad \s_{0}:=\sign(\e),
\ee
where $Q_{k,j} = \FF_{k,j}^{(0)}\s_{0}^{k}$.

%%%%%%%%%%%%%%%%%%%%%%%%%%%%%%%%%%%%%%%%%%%%%%%%%%%%%%%%%%%%%%%%%%%%%%%%%
\begin{lemma}\label{lem:4}
With the notation introduced before, let $\P_{i}^{(0)}$ be
the projection of the segment $\PP_{i}^{(0)}$ on the $j$-axis
and let $\ell_{i}=\ell(\P_{i}^{(0)})$ be the length of $\P_{i}^{(0)}$.
Then $P_{i}^{(0)}(c)$ has $\ell_{i}$ complex non-zero
roots counting multiplicity.
\end{lemma}
%%%%%%%%%%%%%%%%%%%%%%%%%%%%%%%%%%%%%%%%%%%%%%%%%%%%%%%%%%%%%%%%%%%%%%%%%

%%%%%%%%%%%%%%%%%%%%%%%%%%%%%%%%%%%%%%%%%%%%%%%%%%%%%%%%%%%%%%%%%%%%%%%%%
\prova Let $m,n$ be respectively the maximum and the minimum among
the exponents of the variable $\b_{0}$ in $\tilde{\FF}_{i}^{(0)}$.
Then $\ell_{i}=m-n$. Hence $P_{i}^{(0)}$ is a polynomial of
degree $m$ and minimum power $n$: we can write $P_{i}^{(0)}(c) =
c^{n} \tilde{P}(c)$, where $\tilde{P}$ has degree $\ell_{i}$
and $\tilde{P}(0)\neq0$. Fundamental theorem of algebra guarantees
that $\tilde{P}(c)=0$ has $\ell_{i}$ complex solutions counting
multiplicity, which are all the non-zero roots of $P_{i}^{(0)}$.\qed
%%%%%%%%%%%%%%%%%%%%%%%%%%%%%%%%%%%%%%%%%%%%%%%%%%%%%%%%%%%%%%%%%%%%%%%%%
\salto
%%%%%%%%%%%%%%%%%%%%%%%%%%%%%%%%%%%%%%%%%%%%%%%%%%%%%%%%%%%%%%%%%%%%%%%%%

Let $\Re_{0}$ be the set of all the non-zero real solutions of the
polynomial equations $P^{(0)}_{i}(c)=0$. If $\Re_{0}=\emptyset$ the
system (\ref{eq:2.1}) has no subharmonic solution,
as one can easily verify.

Let us suppose then that there exists $c_0\in\Re_{0}$,
such that $c_0(\s_{0}\e)^{\m_{i}^{(0)}}$ is a first
approximate solution of the implicit equation $\FF^{(0)}(\e,\b_{0})=0$
for a suitable $i=1,\ldots,N_{0}$. From now on we shall drop
the label $i$ to lighten the notation. We now set
$\e_{1}=(\s_{0}\e)^{1/\ppp^{(0)}}$, and, as $\e_{1}^{\sss^{(0)}}$
divides $\FF^{(0)}(\s_{0}\e_{1}^{\ppp^{(0)}},c_{0}\e_{1}^{\hhh^{(0)}}+
y_{1}\e_{1}^{\hhh^{(0)}})$, we obtain a new power series
${\eknu \FF 1 \,}(\e_{1},y_{1})$ given by
\be\label{eq:3.6}
\FF^{(0)}(\s_{0}\e_{1}^{\ppp^{(0)}},c_{0}\e_{1}^{\hhh^{(0)}} +
y_{1}\e_{1}^{\hhh^{(0)}}) = \e_{1}^{\sss^{(0)}}{\eknu \FF 1 \,}
(\e_{1},y_{1}) ,
\ee
which is $y_{1}$-general of order $\nnn_{1}$ for some $\nnn_{1}\geq1$.

%%%%%%%%%%%%%%%%%%%%%%%%%%%%%%%%%%%%%%%%%%%%%%%%%%%%%%%%%%%%%%%%%%%%%%%%%
\begin{lemma}\label{lem:5}
With the notations introduced before, let us write
$P^{(0)}(c) = g_{0}(c)(c-c_{0})^{m_{0}}$ with $g_{0}(c_{0}) \neq 0$
and $m_{0} \leq \nnn$. Then $\nnn_{1} = m_{0}$.
\end{lemma}
%%%%%%%%%%%%%%%%%%%%%%%%%%%%%%%%%%%%%%%%%%%%%%%%%%%%%%%%%%%%%%%%%%%%%%%%%

%%%%%%%%%%%%%%%%%%%%%%%%%%%%%%%%%%%%%%%%%%%%%%%%%%%%%%%%%%%%%%%%%%%%%%%%%
\prova This simply follows by the definitions of ${\FF}^{(1)}$
and $P^{(0)}$. In fact we have
\be\label{eq:3.7}
\e_{1}^{\sss^{(0)}}\FF^{(1)}(\e_{1},y_{1})=
\e_{1}^{\sss^{(0)}}\left(\sum_{k+\m^{(0)}j=
\rrr^{(0)}}
\!\!\!\!\!\!
Q_{k,j}(c_{0}+y_{1})^{j}+\e_{1}(\ldots)\right),
\ee
so that $\FF^{(1)}(0,y_{1})=P^{(0)}(c_{0}+y_{1})=g_{0}(c_{0}+
y_{1})y_{1}^{m_{0}}$, and $g_{0}(c_{0}+y_{1})\neq0$ for $y_{1}=0$.
Hence $\FF^{(1)}$ is $y_{1}$-general of order $\nnn_{1}=m_{0}$.\qed
%%%%%%%%%%%%%%%%%%%%%%%%%%%%%%%%%%%%%%%%%%%%%%%%%%%%%%%%%%%%%%%%%%%%%%%%%
\salto
%%%%%%%%%%%%%%%%%%%%%%%%%%%%%%%%%%%%%%%%%%%%%%%%%%%%%%%%%%%%%%%%%%%%%%%%%

Now we restart the process just described: we construct the
Newton polygon $\PP^{(1)}$ of $\FF^{(1)}$. If $\FF^{(1)}_{k,0}=0$
for all $k\geq0$, then $\FF^{(1)}(\e_{1},0)\equiv0$, so that
we have $\FF^{(0)}(\e,c_{0}(\s_{0}\e)^{\m^{(0)}})\equiv0$, \ie
$c_{0}(\s_{0}\e)^{\m^{(0)}}$ is a solution
of the implicit equation $\FF^{(0)}(\e,\b_{0})=0$. Otherwise
we consider the segments $\PP_{1}^{(1)},\ldots,\PP_{N_{1}}^{(1)}$
with slopes $-1/\m_{i}^{(1)}$ for all $i=1,\ldots,N_{1}$, and
we obtain
\be\label{eq:3.8}
\FF^{(1)}(\e_1,y_1)=\tilde{\FF}_{i}^{(1)}(\e_1,y_1) +
\GG_{i}^{(1)}(\e_1,y_1) =
\!\!\!\!\!\!
\sum_{k+j\m_{i}^{(1)} = \rrr_{i}^{(1)}}
\!\!\!\!\!\!
\FF_{k,j}^{(1)}\e_1^{k}y_1^{j} \,\, +
\!\!\!\!\!\!
\sum_{k+j\m_{i}^{(1)} > \rrr_{i}^{(1)}}
\!\!\!\!\!\!
\FF_{k,j}^{(1)}\e_1^{k}y_1^{j},
\ee
where $\rrr_{i}^{(1)}$ is the intercept on the $k$-axis of the
continuation of $\PP_i^{(1)}$.
Hence the first approximate solutions of $\FF^{(1)}(\e_1,y_{1})=0$
are the solutions of the quasi-homogeneous equations
\be\label{eq:3.9}
\tilde{\FF}_{i}^{(1)} (\e_1,y_1) =
\!\!\!\!\!\!
\sum_{k\ppp_{i}^{(1)}+j\hhh_{i}^{(1)}= \sss_{i}^{(1)}}
\!\!\!\!\!\!
\FF_{k,j}^{(1)}\e_1^{k}y_{1}^{j}=0, \qquad i=1,\ldots,N_{1},
\ee
where $\hhh_{i}^{(1)}/\ppp_{i}^{(1)}=\m_{i}^{(1)}$, with
$\hhh_{i}^{(1)}, \ppp_{i}^{(1)}$ relatively prime integers,
and $\sss_{i}^{(1)} = \ppp_{i}^{(1)} \rrr_{i}^{(1)}$.

Thus we define
the polynomials ${\eknu P 1 i}$ such that
\be\label{eq:3.10}
\tilde{\FF}_{i}^{(1)} (\e_1,c\,\e_1^{\m_{i}^{(1)}}) =
\e_1^{\rrr_{i}^{(1)}}
\!\!\!\!\!\!
\sum_{k\ppp_{i}^{(1)}+j\hhh_{i}^{(1)}=\sss_{i}^{(1)}}
\!\!\!\!\!\!
\FF_{k,j}^{(1)}c^{j}=\e_1^{\rrr_{i}^{(1)}}
{\eknu P 1 i}(c),
\ee
and
we call $\Re_{1}$ the set of
the real roots of the polynomials ${\eknu P 1 i}$. If $\Re_{1}=
\emptyset$, we stop the process as there is no subharmonic solution.
Otherwise we call $P^{(1)}$
(\ie again we omit the label $i$) the polynomial which has
a real root $c_1$,
so that $c_1\e_{1}^{\m^{(1)}}$ is an approximate solution
of the equation $\FF^{(1)}(\e_{1},y_{1})=0$. Again
we substitute $\e_{2} = \e_{1}^{1/\ppp^{(1)}}$, and we obtain
\be\label{eq:3.11}
\FF^{(1)}(\e_{2}^{\ppp^{(1)}},c_{1}\e_{2}^{\hhh^{(1)}}+y_{2}
\e_{2}^{\hhh^{(1)}})=
\e_{2}^{\sss^{(1)}}{\eknu \FF 2 \,}(\e_{2},y_{2}) ,
\ee
which is $y_{2}$-general of order $\nnn_{2}\leq\nnn_{1}$, and so on.
Iterating the process we eventually
obtain a sequence of approximate solutions
\be\label{eq:3.12}
\b_{0} =(\s_0\e)^{\m^{(0)}}(c_{0}+y_{1}), \qquad
y_{1} = \e_{1}^{\m^{(1)}}(c_{1}+y_{2}), \qquad
y_{2} = \e_{2}^{\m^{(2)}}(c_{2}+y_3), \qquad \ldots
\ee
where $c_n$ is a (real) root of the polynomial $P^{(n)}$ such that
\be
\FF^{(n)}(\e_n,c\,\e_n^{\m^{(n)}})=\e_n^{\rrr^{(n)}}\left(
P^{(n)}(c)+\e_n(\ldots)
\right),
\label{eq:3.13}\ee
for all $n\ge 0$, where the functions $\FF^{(n)}(\e_{n},y_{n})$
are defined recursively as $\FF^{(n)}(\e^{\ppp^{(n)}}_{n+1},
c_{n}\e^{\hhh^{(n)}}_{n+1}+y_{n+1}\e^{\hhh^{(n)}}_{n+1})=
\e^{\sss^{(n)})}_{n+1}\FF^{(n+1)}(\e_{n+1},y_{n+1})$ for $n\geq1$,
with $\e_{n+1}=\e_{n}^{1/\ppp^{(n)}}$ and the constants
$\m^{(n)}$, $\rrr^{(n)}$, $\sss^{(n)}$, $\hhh^{(n)}$, $\ppp^{(n)}$
defined as in the case $n=0$ in terms of a segment $\PP^{(n)}_{i}$
of the Newton polygon of $\FF^{(n)}$. Therefore
\be
\b_{0} = c_{0}(\s_0\e)^{\m^{(0)}} + c_{1}(\s_0\e)^{\m^{(0)} +
\m^{(1)}/\ppp^{(0)}} + c_{2}(\s_0\e)^{\m^{(0)}+\m^{(1)}/\ppp^{(0)} +
\m^{(2)}/\ppp^{(0)}\ppp^{(1)}} + \ldots
\label{eq:3.14} \ee
is a formal expansion of $\b_{0}$ as a series in ascending fractional
powers of $\s_0\e$.
This iterating method is called the \emph{Newton-Puiseux process}.
Of course this does not occur if we have $\Re_{n}=\emptyset$ at a
certain step ${n}$-th, with ${n}\geq0$.

 From now on we shall suppose $\Re_n\neq\emptyset$ for all $n\geq 0$.
Set also $\nnn_{0}=\nnn$.

%%%%%%%%%%%%%%%%%%%%%%%%%%%%%%%%%%%%%%%%%%%%%%%%%%%%%%%%%%%%%%%%%%%%%%%%%
\begin{lemma}\label{lem:6}
With the notation introduced before, if $\nnn_{i+1}=d_{i}:=\deg(P^{(i)})$
for some $i$, then $\m^{(i)}$ is integer.
\end{lemma}
%%%%%%%%%%%%%%%%%%%%%%%%%%%%%%%%%%%%%%%%%%%%%%%%%%%%%%%%%%%%%%%%%%%%%%%%%

%%%%%%%%%%%%%%%%%%%%%%%%%%%%%%%%%%%%%%%%%%%%%%%%%%%%%%%%%%%%%%%%%%%%%%%%%
\prova Without loss of generality we shall prove the result for
the case $i=0$. Recall that
\be
\FF^{(1)}(0,y_{1})=\sum_{k+\m^{(0)}j=\rrr^{(0)}}Q_{k,j}(c_{0}+y_{1})^{j}
=P^{(0)}(c_{0}+y_{1}),
\label{eq:3.15} \ee
with $\rrr^{(0)}=\m^{(0)}d_{0}$.
If $d_{0}=\nnn_{1}$, then $P^{(0)}$ is of the form
$P^{(0)}(c)=R_{0}(c-c_{0})^{\nnn_{1}}$, with $R_{0}\neq0$.
In particular this means that $Q_{k,\nnn_{1}-1}\neq0$ for some
integer $k\geq0$ with the constraint $k+\m^{(0)}(\nnn_{1}-1)=
\m^{(0)}\nnn_{1}$. Hence $\m^{(0)}=k$ is integer.\qed
%%%%%%%%%%%%%%%%%%%%%%%%%%%%%%%%%%%%%%%%%%%%%%%%%%%%%%%%%%%%%%%%%%%%%%%%%

%%%%%%%%%%%%%%%%%%%%%%%%%%%%%%%%%%%%%%%%%%%%%%%%%%%%%%%%%%%%%%%%%%%%%%%%%
\begin{lemma}\label{lem:7}
With the notations introduced before, there exists $i_{0}\geq0$ such that
$\m^{(i)}$ is integer for all $i\geq i_{0}$.
\end{lemma}
%%%%%%%%%%%%%%%%%%%%%%%%%%%%%%%%%%%%%%%%%%%%%%%%%%%%%%%%%%%%%%%%%%%%%%%%%

%%%%%%%%%%%%%%%%%%%%%%%%%%%%%%%%%%%%%%%%%%%%%%%%%%%%%%%%%%%%%%%%%%%%%%%%%
\prova The series $\FF^{(i)}$ are $y_{i}$-general of order $\nnn_{i}$,
and the $\nnn_{i}$ and the $d_{i}$ form
a descending sequence of natural numbers
\be\label{eq:3.16}
\nnn=\nnn_{0}\geq d_{0}\geq\nnn_{1}\geq d_{1}\geq\ldots
\ee
By Lemma \ref{lem:6}, $\m^{(i)}$ fails to be integers only if
$d_{i}>\nnn_{i+1}$, and this may happen only finitely often.
Hence from a certain $i_{0}$ onwards all the $\m^{(i)}$ are integers.\qed
%%%%%%%%%%%%%%%%%%%%%%%%%%%%%%%%%%%%%%%%%%%%%%%%%%%%%%%%%%%%%%%%%%%%%%%%%
\salto
%%%%%%%%%%%%%%%%%%%%%%%%%%%%%%%%%%%%%%%%%%%%%%%%%%%%%%%%%%%%%%%%%%%%%%%%%

By the results above, we can define $\ppp:=\ppp^{(0)}\cdot\ldots
\cdot\ppp^{(i_{0})}$ such that we can write (\ref{eq:3.14}) as
\be\label{eq:3.17}
\b_{0}=\b_{0}(\e)=\sum_{h\geq h_{0}}{\hknu \b h 0}(\s_0\e)^{h/\ppp},
\ee
where $h_{0}=\hhh^{(0)}\ppp^{(1)}\cdot\ldots\cdot\ppp^{(i_{0})}$.
By construction $\FF^{(0)}(\e,\b_{0}(\e))$ vanishes to all orders,
so that (\ref{eq:3.17}) is a formal solution of
the implicit equation $\FF^{(0)}(\e,\b_{0})=0$.
We shall say that (\ref{eq:3.17}) is a \emph{Puiseux series}
for the plane algebroid curve defined by $\FF^{(0)}(\e,\b_{0})=0$.

%%%%%%%%%%%%%%%%%%%%%%%%%%%%%%%%%%%%%%%%%%%%%%%%%%%%%%%%%%%%%%%%%%%%%%%%%
\begin{lemma}\label{lem:8}
For all $i\geq0$ we can bound $\ppp^{(i)}\leq\nnn_{i}$.
\end{lemma}
%%%%%%%%%%%%%%%%%%%%%%%%%%%%%%%%%%%%%%%%%%%%%%%%%%%%%%%%%%%%%%%%%%%%%%%%%

%%%%%%%%%%%%%%%%%%%%%%%%%%%%%%%%%%%%%%%%%%%%%%%%%%%%%%%%%%%%%%%%%%%%%%%%%
\prova Without loss of generality we prove the result for $i=0$.
By definition, there exist $k',j'$ integers, with $j'\leq\nnn_{0}$,
such that
\be\label{eq:3.18}
\frac{\hhh^{(0)}}{\ppp^{(0)}}=\m^{(0)}=
\frac{\rrr^{(0)}-k'}{j'},
\ee
and $\hhh^{(0)},\ppp^{(0)}$ are relatively prime integers, so that
$\ppp^{(0)}\leq j'\leq\nnn_{0}$.\qed
%%%%%%%%%%%%%%%%%%%%%%%%%%%%%%%%%%%%%%%%%%%%%%%%%%%%%%%%%%%%%%%%%%%%%%%%%
\salto
%%%%%%%%%%%%%%%%%%%%%%%%%%%%%%%%%%%%%%%%%%%%%%%%%%%%%%%%%%%%%%%%%%%%%%%%%

Note that by Lemma \ref{lem:8} we can bound
$\ppp\leq\nnn_{0}\cdot\ldots\cdot\nnn_{i_{0}}\leq\nnn_{0}!$.

%%%%%%%%%%%%%%%%%%%%%%%%%%%%%%%%%%%%%%%%%%%%%%%%%%%%%%%%%%%%%%%%%%%%%%%%%
\begin{lemma}\label{lem:9}
Let $\FF^{(0)}(\e,\b_{0})\in\RRR\{\e,\b_{0}\}$ be $\b_{0}$-general of
order $\nnn$ and let us suppose that $\Re_n\neq\emptyset$ for all
$n\geq0$. Then the series (\ref{eq:3.17}), which formally solves
$\FF^{(0)}(\e,\b_{0}(\e))\equiv0$, is convergent for $\e$ small enough.
\end{lemma}
%%%%%%%%%%%%%%%%%%%%%%%%%%%%%%%%%%%%%%%%%%%%%%%%%%%%%%%%%%%%%%%%%%%%%%%%%

%%%%%%%%%%%%%%%%%%%%%%%%%%%%%%%%%%%%%%%%%%%%%%%%%%%%%%%%%%%%%%%%%%%%%%%%%
\prova Let $P_{\FF^{(0)}}(\e;\b_{0})$ be the Weierstrass polynomial
\cite{BK} of $\FF^{(0)}$ in $\CCC\{\e\}[\b_{0}]$. If $\FF^{(0)}$
is irreducible in $\CCC\{\e,\b_{0}\}$, then by Theorem 1, p. 386
in \cite{BK}, we have a convergent series $\b_{0}(\e^{1/\nnn})$ which
solves the equation $\FF^{(0)}(\e,\b_{0})=0$. Then all the following
\be
\b_{0}\left(\e^{1/\nnn}\right),\;\b_{0}
\left((e^{2\p i}\e)^{1/\nnn}\right),\ldots,
\b_{0}\left((e^{2\p(\nnn-1)i}\e)^{1/\nnn}\right)
\label{eq:3.19} \ee
are solutions of the equation $\FF^{(0)}(\e,\b_{0})=0$.
Thus we have $\nnn$ distinct roots of the Weierstrass polynomial
$P_{\FF^{(0)}}$ and they are all convergent series in
$\CCC\{\e^{1/\nnn}\}$. But also the series (\ref{eq:3.17}) is
a solution of the equation $\FF^{(0)}(\e,\b_{0})=0$. Then, as a polynomial
of degree $\nnn$ has
exactly $\nnn$ (complex) roots counting multiplicity,
(\ref{eq:3.17}) is one of the (\ref{eq:3.19}).
In particular this means that (\ref{eq:3.17}) is convergent
for $\e$ small enough.

In general, we can write
\be\label{eq:3.20}
\FF^{(0)}(\e,\b_{0})=\prod_{i=1}^N(\FF^{(0)}_{i}(\e,\b_{0}))^{m_{i}},
\ee
for some $N\geq1$, where the $\FF^{(0)}_{i}$ are the irreducible factors
of $\FF^{(0)}$. Then the Puiseux series (\ref{eq:3.17}) solves
one of the equations $\FF^{(0)}_{i}(\e,\b_{0})=0$, and hence,
by what said above, it converges for $\e$ small enough.\qed
%%%%%%%%%%%%%%%%%%%%%%%%%%%%%%%%%%%%%%%%%%%%%%%%%%%%%%%%%%%%%%%%%%%%%%%%%
\salto
%%%%%%%%%%%%%%%%%%%%%%%%%%%%%%%%%%%%%%%%%%%%%%%%%%%%%%%%%%%%%%%%%%%%%%%%%

As a consequence of Lemma \ref{lem:9} we obtain the following corollary.

%%%%%%%%%%%%%%%%%%%%%%%%%%%%%%%%%%%%%%%%%%%%%%%%%%%%%%%%%%%%%%%%%%%%%%%%%
\begin{teo}\label{thm:1}
Consider a periodic solution with frequency $\o=p/q$ for the system
(\ref{eq:2.1}). Assume that
Hypotheses \ref{hyp1} and \ref{hyp2} are satisfied with $\nnn$ odd.
Then for $\e$ small enough the system (\ref{eq:2.1})
has at least one subharmonic solution of order $q/p$.
Such a solution admits a convergent power series in $|\e|^{1/\nnn!}$,
and hence a convergent Puiseux series in $|\e|$.
\end{teo}
%%%%%%%%%%%%%%%%%%%%%%%%%%%%%%%%%%%%%%%%%%%%%%%%%%%%%%%%%%%%%%%%%%%%%%%%%

%%%%%%%%%%%%%%%%%%%%%%%%%%%%%%%%%%%%%%%%%%%%%%%%%%%%%%%%%%%%%%%%%%%%%%%%%
\prova If $\nnn$ is odd, then $\Re_n\neq\emptyset$ for all $n\geq0$.
This trivially follows from the fact that if $\nnn$ is odd,
then there exists at least one polynomial $P_{i}^{(0)}$ associated
with a segment $\PP_{i}^{(0)}$ whose projection $\P_{i}^{(0)}$
on the $j$-axis is associated with a polynomial $\tilde{P}_{i}^{(0)}$
with odd degree $\ell_{i}$. Thus such a polynomial admits a non-zero
real root with odd multiplicity $\nnn_{1}$, so that
$\FF^{(1)}(\e_{1},y_{1})$
is $y_{1}$-general of odd order $\nnn_{1}$ and so on.

Hence we can apply the Newton-Puiseux process to obtain a subharmonic
solution as a Puiseux series in $\e$ which is convergent for $\e$
sufficiently small by Lemma \ref{lem:9}.\qed
%%%%%%%%%%%%%%%%%%%%%%%%%%%%%%%%%%%%%%%%%%%%%%%%%%%%%%%%%%%%%%%%%%%%%%%%%
\salto
%%%%%%%%%%%%%%%%%%%%%%%%%%%%%%%%%%%%%%%%%%%%%%%%%%%%%%%%%%%%%%%%%%%%%%%%%

Theorem \ref{thm:1} extends the results of \cite{ZL}.
First it gives the explicit dependence of the parameter $\b_{0}$ on $\e$,
showing that it is analytic in $|\e|^{1/\nnn!}$.
Second, it shows that it is possible to express
the subharmonic solution as a convergent fractional power series
in $\e$, and this allows us to push perturbation theory
to arbitrarily high order. 

On the other hand, the Newton-Puiseux algorithm does not allow
to construct the solution within any fixed accuracy. In this regard,
it is not really constructive: we know that the solution is analytic
in a fractional power of $\e$, but we know neither the
size of the radius of convergence nor the precision with
which the solution is approximated if we stop
the Newton-Puiseux at a given step. Moreover we know that there is
at least one subharmonic solution, but we are not able to
decide how many of them are possible. In fact, a subharmonic solution
can be constructed for any non-zero real root of each odd-degree
polynomial $P_{i}^{(n)}$ associated with each segment of $\PP^{(n)}$
to all step of iteration, but we cannot predict \emph{a priori}
how many possibilities will arise along the process.

However, we obtain a fully constructive algorithm if
we make some further hypothesis.

%%%%%%%%%%%%%%%%%%%%%%%%%%%%%%%%%%%%%%%%%%%%%%%%%%%%%%%%%%%%%%%%%%%%%%%%%
\begin{hyp}\label{hyp3}
There exists $i_{0}\geq0$ such that at the $i_{0}$-th step of the iteration,
there exists a polynomial
$P^{(i_{0})}=P^{(i_{0})}(c)$ which has a
simple root $c_{i_0}\in\RRR$.
\end{hyp}
%%%%%%%%%%%%%%%%%%%%%%%%%%%%%%%%%%%%%%%%%%%%%%%%%%%%%%%%%%%%%%%%%%%%%%%%%

Indeed, if we assume Hypothesis \ref{hyp3},
we obtain the following result.

%%%%%%%%%%%%%%%%%%%%%%%%%%%%%%%%%%%%%%%%%%%%%%%%%%%%%%%%%%%%%%%%%%%%%%%%%
\begin{teo}\label{thm:2}
Consider a periodic solution with frequency $\o=p/q$ for the system
(\ref{eq:2.1}). Assume that Hypotheses \ref{hyp1}, \ref{hyp2}
and \ref{hyp3} are satisfied. Then there exists an
explicitly computable value $\e_{0}>0$ such that for $|\e|<\e_{0}$
the system (\ref{eq:2.1}) has at least
one subharmonic solution of order $q/p$.
Such a solution admits a convergent power series in $|\e|^{1/\nnn!}$,
and hence a convergent Puiseux series in $|\e|$.
\end{teo}
%%%%%%%%%%%%%%%%%%%%%%%%%%%%%%%%%%%%%%%%%%%%%%%%%%%%%%%%%%%%%%%%%%%%%%%%%

We shall see in Section \ref{sec:4} that, by assuming
Hypothesis \ref{hyp3}, we can use the Newton-Puiseux algorithm
up to the $i_{0}$-th step (hence a finite number of times), and we
can provide recursive formulae for the higher order contributions.
This will allow us -- as we shall see in Section \ref{sec:6} -- to
introduce a graphical representation for the subharmonic solution,
and, eventually, to obtain an explicit bound
on the radius of convergence of the power series expansion.

%%%%%%%%%%%%%%%%%%%%%%%%%%%%%%%%%%%%%%%%%%%%%%%%%%%%%%%%%%%%%%%%%%%%%%%%%
%%%%%%%%%%%%%%%%%%%%%%%%%%%%%%%%%%%%%%%%%%%%%%%%%%%%%%%%%%%%%%%%%%%%%%%%%
\zerarcounters
\section{Trees expansion and proof of Lemma \ref{lem:1}}
\label{sec:5}
%%%%%%%%%%%%%%%%%%%%%%%%%%%%%%%%%%%%%%%%%%%%%%%%%%%%%%%%%%%%%%%%%%%%%%%%%
%%%%%%%%%%%%%%%%%%%%%%%%%%%%%%%%%%%%%%%%%%%%%%%%%%%%%%%%%%%%%%%%%%%%%%%%%

\noindent A tree $\th$ is defined as a partially ordered set
of points $\vvv$ (\emph{vertices}) connected by
oriented \emph{lines} $\ell$. The lines are consistently
oriented toward a unique point called the \emph{root} which
admits only one entering line called the \emph{root line}.
If a line $\ell$ connects two vertices $\vvv_{1},\vvv_{2}$ and is oriented
from $\vvv_{2}$ to $\vvv_{1}$, we say that
$\vvv_{2}\prec\vvv_{1}$ and we shall write $\ell_{\vvv_{2}}=\ell$.
We shall say that $\ell$ exits $\vvv_{2}$ and enters $\vvv_{1}$.
More generally we write $\vvv_{2}\prec\vvv_{1}$ when $\vvv_{1}$ is
on the path of lines
connecting $\vvv_{2}$ to the root: hence the orientation of the
lines is opposite to the partial ordering
relation $\prec$.

We denote with $V(\th)$ and $L(\th)$ the set of vertices and
lines in $\th$ respectively,
and with $|V(\th)|$ and $|L(\th)|$ the number of vertices and
lines respectively.
Remark that one has $|V(\th)|=|L(\th)|$.

We consider two kinds of vertices: \emph{nodes} and \emph{leaves}.
The leaves can only be end-points, \ie points with no lines
entering them, while the nodes can be either end-points or not.
We shall not consider the tree consisting of only one leaf and
the line exiting it, \ie a tree must have at least the node
which the root line exits.

We shall denote with $N(\th)$ and $E(\th)$ the set of nodes
and leaves respectively. Here and henceforth we shall denote
with $\vvv$ and $\eee$ the nodes and the leaves respectively.
Remark that $V(\th)=N(\th)\amalg E(\th)$.

With each line $\ell=\ell_\vvv$, we associate three labels
$(h_\ell,\d_\ell,\n_\ell)$, with $h_\ell\in\{\a,A\}$,
$\d_{\ell}\in\{1,2\}$ and $\n_\ell\in\ZZZ$, with the constraint
that $\n_\ell\neq0$ for $h_\ell=\a$ and $\d_\ell=1$ for $h_\ell=A$.
With each line $\ell=\ell_\eee$ we associate $h_\ell=\a$,
$\d_\ell=1$ and $\n_\ell=0$.
We shall say that $h_{\ell}$, $\d_{\ell}$ and $\n_\ell$
are the \emph{component} label, the \emph{degree} label
and the \emph{momentum} of the line $\ell$, respectively.

Given a node $\vvv$, we call $r_{\vvv}$ the number
of the lines entering $\vvv$ carrying a component label $h=\a$
and $s_{\vvv}$ the number of the lines entering $\vvv$
with component label $h=A$. We also introduce a \emph{badge}
label $b_\vvv\in\{0,1\}$ with the constraint that
$b_{\vvv}=1$ for $h_{\ell_{\vvv}}=\a$ and $\d_{\ell_{\vvv}}=2$,
and for $h_{\ell_\vvv}=A$ and $\n_{\ell_\vvv}\neq0$,
and two \emph{mode} labels
$\s_{\vvv},\s'_{\vvv}\in\ZZZ$. 
We call \emph{global mode} label the sum
\be\label{eq:5.1}
\n_{\vvv}=p\s_{\vvv}+q\s'_{\vvv},
\ee
where $q,p$ are the relatively prime integers such that $\o(A_{0})=p/q$,
with the constraint that $\n_{\vvv}=0$ when
$b_{\vvv}=0.$

For all $\ell=\ell_\vvv$,
we set also the following \emph{conservation law}
\be\label{eq:5.2}
\n_{\ell}=\n_{\ell_{\vvv}}=\sum_{\substack{\www\,\in\,N(\th) \\
\www\,\preceq\,\vvv}}\n_{\www},
\ee
\ie the momentum of the line exiting $\vvv$ is
the sum of the momenta of the lines entering $\vvv$ plus the global
mode of the node $\vvv$ itself.

Given a labeled tree $\th$, where labels are defined as above,
we associate with each line $\ell$ exiting a node,
a \emph{propagator}
\be\label{eq:5.3}
g_{\ell}=\left\{\ba
	& \frac{\o'(A_{0})^{\d_{\ell}-1}}{(i\o\n_{\ell})^{\d_{\ell}}},
&h_{\ell}=\a,A,\quad \n_{\ell}\neq0, \\
	&-\frac{1}{\o'(A_{0})}, &h_{\ell}=A,\quad\quad \n_{\ell}=0,
\ea\right.\ee
while for each line $\ell$ exiting a leaf we set $g_{\ell}=1$.

Moreover, we associate with each node $\vvv$ a \emph{node factor}
\be\label{eq:5.4}
\NN_{\vvv}=\left\{\ba
& \frac{(i\s_{\vvv})^{r_{\vvv}}\dpr_A^{s_{\vvv}}}{r_{\vvv}!s_{\vvv}!}
F_{\s_{\vvv},\s'_{\vvv}}(A_{0},t_{0}),
	&h_{\ell_\vvv}=\a,\quad\d_{\ell_{\vvv}}=1,\quad b_{\vvv}=1,
\quad \n_{\ell_\vvv}\neq0,\\
& \frac{\dpr_A^{s_{\vvv}}}{s_{\vvv}!}\o(A_{0}), 
	&h_{\ell_\vvv}=\a,\quad\d_{\ell_{\vvv}}=1,\quad b_{\vvv}=0,
\quad\n_{\ell_\vvv}\neq0,\\
& \frac{(i\s_{\vvv})^{r_{\vvv}}\dpr_A^{s_{\vvv}}}{r_{\vvv}!s_{\vvv}!}
G_{\s_{\vvv},\s'_{\vvv}}(A_{0},t_{0}),
	&h_{\ell_\vvv}=\a,\quad\d_{\ell_{\vvv}}=2,\quad b_{\vvv}=1,
\quad\n_{\ell_\vvv}\neq0,\\
& \frac{(i\s_{\vvv})^{r_{\vvv}}\dpr_A^{s_{\vvv}}}{r_{\vvv}!s_{\vvv}!}
G_{\s_{\vvv},\s'_{\vvv}}(A_{0},t_{0}),
	&h_{\ell_\vvv}=A,\quad\d_{\ell_{\vvv}}=1,\quad b_{\vvv}=1,
\quad\n_{\ell_\vvv}\neq0,\\
& \frac{(i\s_{\vvv})^{r_{\vvv}}\dpr_A^{s_{\vvv}}}{r_{\vvv}!s_{\vvv}!}
F_{\s_{\vvv},\s'_{\vvv}}(A_{0},t_{0}),
	&h_{\ell_\vvv}=A,\quad\d_{\ell_{\vvv}}=1,\quad b_{\vvv}=1,
\quad \n_{\ell_\vvv}=0,\\
& \frac{\dpr_A^{s_{\vvv}}}{s_{\vvv}!}\o(A_{0}), 
	&h_{\ell_\vvv}=A,\quad\d_{\ell_{\vvv}}=1,\quad b_{\vvv}=0,
\quad\n_{\ell_\vvv}=0,
\ea\right.\ee
with the constraint that when $b_{\vvv}=0$ one has $r_{\vvv}=0$
and $s_{\vvv}\geq2$.

Given a labeled tree $\th$ with propagators and node factors
associated as above, we define the
\emph{value} of $\th$ the number
\be\label{eq:5.5}
\Val(\th)=\left(\prod_{\ell\in L(\th)}g_{\ell}\right)
\left(\prod_{\vvv\in N(\th)}\NN_{\vvv}\right).
\ee
Remark that $\Val(\th)$ is a well-defined quantity because all the
propagators and node factors are bounded
quantities.

For each line $\ell$ exiting a node $\vvv$ we set
$b_{\ell}=b_{\vvv}$, while for each line
$\ell$ exiting a leaf we set $b_{\ell}=0$.
Given a labeled tree $\th$, we call \emph{order} of $\th$ the number
\be\label{eq:5.6}
k(\th)=|\{\ell\in L(\th)\;:\;b_{\ell}=1\}|;
\ee
the momentum $\n(\th)$ of the root line will be the
\emph{total momentum}, and the component
label $h(\th)$ associated to the root line will be the
\emph{total component label}.
Moreover, we set $j(\th)=|E(\th)|$.

Define $\TT_{k,\n,h,j}$ the set of all the trees $\th$ with
order $k(\th)=k$, total momentum $\n(\th)=\n$,
total component label $h(\th)=h$ and $j(\th)=j$ leaves.

%%%%%%%%%%%%%%%%%%%%%%%%%%%%%%%%%%%%%%%%%%%%%%%%%%%%%%%%%%%%%%%%%%%%%%%%%
\begin{lemma}\label{lem:10}
For any tree $\th$ labeled as before, one has $|L(\th)|=|V(\th)|
\leq2k(\th)+j(\th)-1$.
\end{lemma}
%%%%%%%%%%%%%%%%%%%%%%%%%%%%%%%%%%%%%%%%%%%%%%%%%%%%%%%%%%%%%%%%%%%%%%%%%

%%%%%%%%%%%%%%%%%%%%%%%%%%%%%%%%%%%%%%%%%%%%%%%%%%%%%%%%%%%%%%%%%%%%%%%%%
\prova We prove the bound $|N(\th)|\leq2k(\th)-1$ by induction on $k$.

For $k=1$ the bound is trivially satisfied, as a direct check shows: in
particular, a tree $\th$ with $k(\th)=1$ has exactly one node and $j(\th)$
leaves.
In fact if $\th$ has a line $\ell=\ell_{\vvv}$ with $b_{\ell}=0$,
then $\vvv$ has $s_{\vvv}\geq2$
lines with component label $h=A$ entering it. Hence there are at least
two lines exiting a node with $b_{\vvv}=1$.

Assume now that the bound holds for all $k'<k$, and let us show
that then it holds also for $k$.
Let $\ell_{0}$ be the root line of $\th$ and $\vvv_{0}$ the node
which the root line exits.
Call $r$ and $s$ the number of lines entering $\vvv_{0}$ with component
labels $\a$ and $A$ respectively,
and denote with $\th_{1},\ldots,\th_{r+s}$ the subtrees which have those
lines as root lines. Then
\be\label{eq:5.7}
|N(\th)|=1+\sum_{m=1}^{r+s}|N(\th_{m})|.
\ee

If $\ell_{0}$ has badge label $b_{\ell_{0}}=1$ we have
$|N(\th)|\leq1+2(k-1)-(r+s)\leq2k-1$,
by the inductive hypothesis and by the fact that $k(\th_{1})+ \ldots+
k(\th_{r+s})=k-1$. If $\ell_{0}$ has badge label $b_{\ell_{0}}=0$ we have
$|N(\th)|\leq1+2k-(r+s)\leq2k-1$,
by the inductive hypothesis, by the fact that $k(\th_{1})+\ldots+
k(\th_{r+s})=k$, and the constraint that $s\geq2$.
Therefore the assertion is proved.\qed
%%%%%%%%%%%%%%%%%%%%%%%%%%%%%%%%%%%%%%%%%%%%%%%%%%%%%%%%%%%%%%%%%%%%%%%%%

%%%%%%%%%%%%%%%%%%%%%%%%%%%%%%%%%%%%%%%%%%%%%%%%%%%%%%%%%%%%%%%%%%%%%%%%%
\begin{lemma}\label{lem:11}
The Fourier coefficients ${\eknu {\ol{\b}} {k,j} \n}$, $\n\neq0$, and
${\eknu {\ol{B}} {k,j} \n}$
can be written in terms of trees as
\begin{subequations}
\begin{align}
{\eknu {\ol{\b}} {k,j} \n} & =\sum_{\th\in\TT_{k,\n,\a,j}}
\Val(\th),\quad\n\neq0,
\label{eq:5.8a} \\
{\eknu {\ol{B}} {k,j} \n} & =\sum_{\th\in\TT_{k,\n,A,j}}
\Val({\th}),\quad\n\in\ZZZ,
\label{eq:5.8b}
\end{align}
\label{eq:5.8} \end{subequations}
\vskip-.3truecm
\noindent for all $k\geq1$, $j\geq0$.
\end{lemma}
%%%%%%%%%%%%%%%%%%%%%%%%%%%%%%%%%%%%%%%%%%%%%%%%%%%%%%%%%%%%%%%%%%%%%%%%%

\begin{figure}[!ht]
\centering
\begin{minipage}[b]{8cm}
%\centering
{
\psfrag{be}{$\ol{\b}_\n^{(k,0)}$}
\psfrag{B}{$\ol{B}_\n^{(k,0)}$}
\psfrag{a}{$\a$}
\psfrag{d}{$\d$}
\psfrag{n}{$\n$}
\psfrag{+}{$+$}
\psfrag{2}{$2$}
\psfrag{1}{$1$}
\psfrag{(k,0)}{$(k,0)$}
\psfrag{A}{$A$}
\includegraphics[width=11cm]{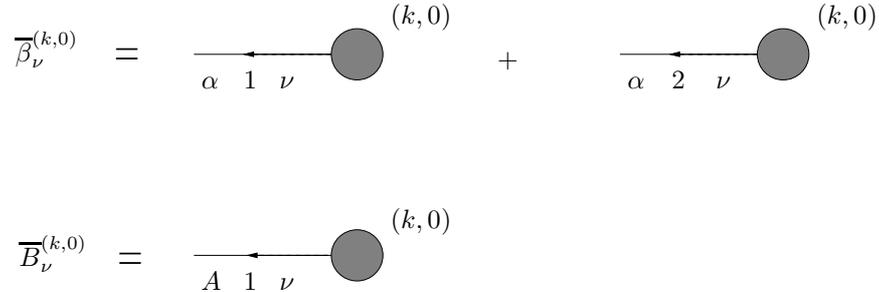}
}
\caption{\label{fig:2}\footnotesize{Graph elements.}}
\end{minipage}
\end{figure}

%%%%%%%%%%%%%%%%%%%%%%%%%%%%%%%%%%%%%%%%%%%%%%%%%%%%%%%%%%%%%%%%%%%%%%%%%
\prova First we consider trees without leaves, \ie
the coefficients $\ol{\b}^{(k,0)}_{\n}$, $\n\neq0$, and
$\ol{B}^{(k,0)}_{\n}$. For $k=1$ is a direct check.
Now let us suppose that the assertion holds for all $k<\ol{k}$.
Let us write $f_\a=\ol{\b}$, $f_A=\ol{B}$ and represent
the coefficients $f^{(k,0)}_{\n,h}$ with the graph elements
in Figure \ref{fig:2}, as a line with label $\n$ and $h=\a,A$
respectively, exiting a ball with label $(k,0)$.

%%%%%%%%%%%%%%%%%%%%%%%%%%%%%%%%%%%%%%%%%%%%%%%%%%%%%%%%%%%%%%%%%%%%%%%%%
% figure 2
%%%%%%%%%%%%%%%%%%%%%%%%%%%%%%%%%%%%%%%%%%%%%%%%%%%%%%%%%%%%%%%%%%%%%%%%%
%%%%%%%%%%%%%%%%%%%%%%%%%%%%%%%%%%%%%%%%%%%%%%%%%%%%%%%%%%%%%%%%%%%%%%%%%

\begin{figure}[!ht]
\centering
\begin{minipage}[b]{10cm}
%\centering
{
\psfrag{h}{$h$}
\psfrag{A}{$A$}
\psfrag{a}{$\a$}
\psfrag{d}{$\d$}
\psfrag{n}{$\n$}
\psfrag{(k,0)}{$(\ol{k},0)$}
\psfrag{v0}{$\vvv_{0}$}
\psfrag{n0}{$\n_{0}$}
\psfrag{d1}{$\d_{1}$}
\psfrag{dr}{$\d_{r}$}
\psfrag{dr+1}{$\d_{r+1}$}
\psfrag{dm}{$\d_{m}$}
\psfrag{(k1,0)}{$(k_{1},0)$}
\psfrag{(kr,0)}{$(k_{r},0)$}
\psfrag{(kr+1,0)}{$(k_{r+1},0)$}
\psfrag{(km,0)}{$(k_{m},0)$}
\psfrag{n1}{$\n_{1}$}
\psfrag{nr}{$\n_{r}$}
\psfrag{nr+1}{$\n_{r+1}$}
\psfrag{nm}{$\n_{m}$}
\includegraphics[width=10cm]{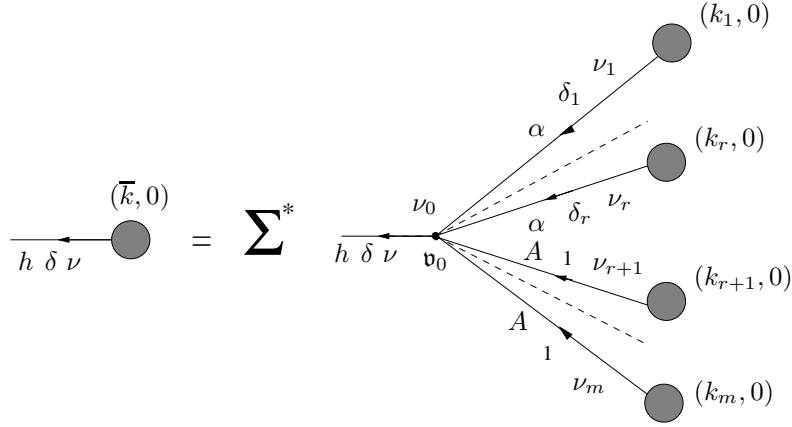}
}
\caption{\label{fig:3}\footnotesize{Graphical representation
for the recursive equations (\ref{eq:2.23}).}}
\end{minipage}
\end{figure}

Then we can represent each equation of (\ref{eq:2.21}) graphically
as in Figure \ref{fig:3}, simply representing each factor
${\eknu f {k_{i},0} {h_{i},\n_{i}}}$ in the r.h.s. as a graph element
according to Figure \ref{fig:2}: the lines of all such graph
elements enter the same node $\vvv_{0}$.

%%%%%%%%%%%%%%%%%%%%%%%%%%%%%%%%%%%%%%%%%%%%%%%%%%%%%%%%%%%%%%%%%%%%%%%%%
% figure 3
%%%%%%%%%%%%%%%%%%%%%%%%%%%%%%%%%%%%%%%%%%%%%%%%%%%%%%%%%%%%%%%%%%%%%%%%%
%%%%%%%%%%%%%%%%%%%%%%%%%%%%%%%%%%%%%%%%%%%%%%%%%%%%%%%%%%%%%%%%%%%%%%%%%

The root line $\ell_{0}$ of such trees will carry a component label
$h=\a,A$ for $f=\ol{\b},\ol{B}$
respectively, and a momentum label $\n$.
Hence, by inductive hypothesis, one obtains
\be\label{eq:5.9}\ba
{\eknu f {\ol{k},0} {h,\n}}&={\mathop\sum}^{*} g_{\ell_{0}}
\NN_{\vvv_{0}}{\eknu f {k_{1},0}{h_{1},\n_{1}}}\ldots
	{\eknu f {k_{m},0}{h_{m},\n_{m}}}\\
&={\mathop\sum}^{*} g_{\ell_{0}}\NN_{\vvv_{0}}
	\left(\sum_{\th\in\TT_{\k_{1},\n_{1},h_{1},0}}
\Val(\th)\right)\ldots
	\left(\sum_{\th\in\TT_{\k_{m},\n_{m},h_{m},0}}
\Val(\th)\right) =\sum_{\th\in\TT_{\ol{k},\n,h,0}}\Val(\th),
\ea\ee
where $m=r_{0}+s_{0}$, and we write $\sum^*$ for the sum over
all the labels admitted by the constraints,
so that the assertion is proved for all $k$ and for $j=0$.

Now we consider $k$ as fixed and we prove the statement by
induction on $j$. The case $j=0$ has
already been discussed.
Finally we assume that the assertion holds for $j=j'$ and 
show that then it holds for $j'+1$.
Notice that a tree $\th\in\TT_{k,\n,h,j'+1}$, for both $h=\a,A$ 
can be obtained by considering
a suitable tree $\th_{0}\in\TT_{k,\n,h,j'}$ attaching an extra 
leaf to a node of $\th_{0}$ and applying an
extra derivative $\dpr_{\a}$ to the node factor associated 
with that node. If one considers all the trees
that can be obtained in such a way from the same $\th_{0}$ 
and sums together all those contributions, one finds
a quantity proportional to $\dpr_{\a}\Val(\th_{0})$. Then if we 
sum over all possible choices of $\th_{0}$,
we reconstruct $\ol{\b}^{(k,j'+1)}_{\n}$ for $h=\a$ 
and $\ol{B}^{(k,j'+1)}_{\n}$ for $h=A$.
Hence the assertion follows.\qed
%%%%%%%%%%%%%%%%%%%%%%%%%%%%%%%%%%%%%%%%%%%%%%%%%%%%%%%%%%%%%%%%%%%%%%%%%

%%%%%%%%%%%%%%%%%%%%%%%%%%%%%%%%%%%%%%%%%%%%%%%%%%%%%%%%%%%%%%%%%%%%%%%%%
\begin{prop}\label{prop:1}
The formal solution (\ref{eq:2.19}) of the system (\ref{eq:2.10}),
given by the recursive
equations (\ref{eq:2.21}),
converges for $\e$ and $\b_{0}$ small enough.
\end{prop}
%%%%%%%%%%%%%%%%%%%%%%%%%%%%%%%%%%%%%%%%%%%%%%%%%%%%%%%%%%%%%%%%%%%%%%%%%

%%%%%%%%%%%%%%%%%%%%%%%%%%%%%%%%%%%%%%%%%%%%%%%%%%%%%%%%%%%%%%%%%%%%%%%%%
\prova First of all we remark that by Lemma \ref{lem:10},
the number of unlabeled trees of order 
$k$ and $j$ leaves is bounded by $4^{2k+j}\times2^{2k+j}=8^{2k+j}$.
The sum over all labels except the mode labels and the momenta is
bounded again by a constant to the power $k$ times a constant to
the power $j$, simply because all such labels can assume only a
finite number of values. Now by the analyticity assumption
on the functions $F$ and $G$, we have
the bound
\be\label{eq:5.10}\ba
&\left| \frac{(i\s_{0})^r}{r!} \frac{\dpr^s_{A}}{s!}
F_{\s_{0},\s'_{0}}(A_{0},t_{0})\right|
\leq \QQ R^rS^se^{-\k(|\s_{0}|+|\s'_{0}|)},\\
&\left| \frac{(i\s_{0})^r}{r!} \frac{\dpr^s_{A}}{s!}
G_{\s_{0},\s'_{0}}(A_{0},t_{0})\right|
\leq \QQ R^rS^se^{-\k(|\s_{0}|+|\s'_{0}|)},
\ea\ee
for suitable positive constants $\QQ,R,S,\k$, and we can imagine,
without loss of generality,
that $\QQ$ and $S$ are such that $|\dpr^s_A\o(A_{0})/s!|\leq \QQ S^s$.
This gives us a bound for the node factors.
The propagators can be bounded by
\be\label{eq:5.11}
|g_{\ell}|\leq\max\left\{\left| \frac{\o'(A_{0})}{\o^2}\right|,
\left|\frac{1}{\o'(A_{0})}\right|,\left|\frac{1}{\o}\right|,1\right\},
\ee
so that the product over all the lines can be bounded again
by a constant to the power $k$ times
a constant to the power $j$.

Thus the sum over the mode labels -- which uniquely determine the
momenta -- can be performed by using for each node half the
exponential decay factor provided by (\ref{eq:5.10}). Then we obtain
\be\label{eq:5.12}
|\ol{\b}^{(k,j)}_{\n}|\leq C_{1}C^{k}_{2}C^{j}_3e^{-\k|\n|/2},\quad
|\ol{B}^{(k,j)}_{\n}|\leq C_{1}C^{k}_{2}C^{j}_3e^{-\k|\n|/2},
\ee
for suitable constants $C_{1}$, $C_{2}$ and $C_3$.
This provides the convergence of the series (\ref{eq:2.19})
for $|\e|<C^{-1}_{2}$ and $|\b_{0}|<C^{-1}_{3}$.\qed
%%%%%%%%%%%%%%%%%%%%%%%%%%%%%%%%%%%%%%%%%%%%%%%%%%%%%%%%%%%%%%%%%%%%%%%%%

%%%%%%%%%%%%%%%%%%%%%%%%%%%%%%%%%%%%%%%%%%%%%%%%%%%%%%%%%%%%%%%%%%%%%%%%%
%%%%%%%%%%%%%%%%%%%%%%%%%%%%%%%%%%%%%%%%%%%%%%%%%%%%%%%%%%%%%%%%%%%%%%%%%
\zerarcounters
\section{Formal solubility of the equations of motion}
\label{sec:4}
%%%%%%%%%%%%%%%%%%%%%%%%%%%%%%%%%%%%%%%%%%%%%%%%%%%%%%%%%%%%%%%%%%%%%%%%%
%%%%%%%%%%%%%%%%%%%%%%%%%%%%%%%%%%%%%%%%%%%%%%%%%%%%%%%%%%%%%%%%%%%%%%%%%

\noindent Assume that Hypotheses \ref{hyp1}, \ref{hyp2} and \ref{hyp3}
are satisfied. Let us set $\h:=|\e|^{1/\ppp}$,
where $\ppp=\ppp^{(0)}\cdot\ldots\cdot\ppp^{(i_{0})}$.
We search for a formal solution $(\a(t),A(t))$ of (\ref{eq:2.1}),
with $\a(t)=\a_{0}(t)+\b_{0}+\tilde{\b}(t)$ and $A(t)=A_{0}+B(t)$, where
\be\label{eq:4.1}
\b_{0}=\sum_{k\geq1}\h^{k}{\hknu \b k 0},\qquad
\bt(t)=\sum_{\substack{\n\in\ZZZ\\ \n\neq0}}e^{i\n\o t}
\sum_{k\geq1}\h^{k}{\hknu {\bt} k \n},\qquad
B(t)=\sum_{\n\in\ZZZ}e^{i\n\o t}\sum_{k\geq1}\h^{k}{\hknu B k \n},
\ee
and the coefficients ${\hknu \b k 0}$, ${\hknu \bt k \n}$
and ${\hknu B k \n}$ solve
\be\label{eq:4.2}
\left\{\ba
\bt^{[k]}_{\n} &=\frac{\F^{[k]}_{\n}}{i\o\n}+
\o'(A_{0}) \frac{\G^{[k]}_{\n}}{(i\o\n)^2},\qquad
B^{[k]}_{\n} =\frac{\G^{[k]}_{\n}}{i\o\n}, \qquad \n\neq0, \\
B^{[k]}_{0} & =-\frac{\F^{[k]}_{0}}{\o'(A_{0})}, \qquad
\G^{[k]}_{0} =0,
\ea\right.\ee
with the functions ${\hknu {\G} k \n}$ and ${\hknu {\F} k \n}$
recursively defined as
\begin{subequations}
\begin{align}
{\G}^{[k]}_{\n} & = \sum_{m\geq0}\,\sum_{r+s=m}\,
\sum_{\substack{p\s_{0}+q\s'_{0}+\n_{1}+\ldots+\n_{m}=\n\\ 
	k_{1}+\ldots+k_{m}=k-\ppp}} \frac{(i\s_{0})^r}{r!}
\frac{\dpr_A^{s}}{s!}G_{\s_{0},\s'_{0}}(A_{0},t_{0})
	{\b}^{[k_{1}]}_{\n_{1}}\cdots{\b}^{[k_{r}]}_{\n_{r}}
{B}^{[k_{r+1}]}_{\n_{r+1}}\cdots
	{B}^{[k_{m}]}_{\n_{m}},
\label{eq:4.3a} \\
{\F}^{[k]}_{\n} & = \sum_{m\geq0}\,\sum_{r+s=m}\,
\sum_{\substack{p\s_{0}+q\s'_{0}+\n_{1}+\ldots+\n_{m}=\n\\ 
	k_{1}+\ldots+k_{m}=k-\ppp}} \frac{(i\s_{0})^r}{r!}
\frac{\dpr_A^{s}}{s!}F_{\s_{0},\s'_{0}}(A_{0},t_{0})
	{\b}^{[k_{1}]}_{\n_{1}}\cdots{\b}^{[k_{r}]}_{\n_{r}}
{B}^{[k_{r+1}]}_{\n_{r+1}}\cdots
	{B}^{[k_{m}]}_{\n_{m}}, \nonumber \\
& +\sum_{s\geq2}\,\sum_{\substack{\n_{1}+\ldots+\n_{s}=\n\\ 
k_{1}+\ldots+k_{s}=k}}
	\frac{\dpr_A^s}{s!}\o(A_{0}){\hknu B {k_{1}} {\n_{1}}}
\ldots{\hknu B {k_{s}} {\n_{s}}},
\label{eq:4.3b}
\end{align}
\label{eq:4.3} \end{subequations}
\vskip-.3truecm
\noindent where ${\hknu \b k \n}={\hknu \bt k \n}$ for $\n\neq0$ .
We use a different notation
for the Taylor coefficients to stress that we are expanding in $\h$.

We say that the integral equations (\ref{eq:2.9}), and
hence the equations (\ref{eq:4.2}), are satisfied up to
order $\ol{k}$ if there exists a choice of the parameters
${\hknu \b 1 0},\ldots{\hknu \b {\ol{k}} 0}$ which make the relations
(\ref{eq:4.2}) to be satisfied for all $k=1,\ldots,\ol{k}$.

%%%%%%%%%%%%%%%%%%%%%%%%%%%%%%%%%%%%%%%%%%%%%%%%%%%%%%%%%%%%%%%%%%%%%%%%%
\begin{lemma}\label{lem:12}
The equations (\ref{eq:4.2}) are satisfied up to order $k=\ppp-1$
with ${\hknu \bt k \n}$ and
${\hknu B k \n}$ identically zero for all $k=1,\ldots,\ppp-1$ and for
any choice of the constants
${\hknu \b 1 0},\ldots{\hknu \b {\ppp-1} 0}$.
\end{lemma}
%%%%%%%%%%%%%%%%%%%%%%%%%%%%%%%%%%%%%%%%%%%%%%%%%%%%%%%%%%%%%%%%%%%%%%%%%

%%%%%%%%%%%%%%%%%%%%%%%%%%%%%%%%%%%%%%%%%%%%%%%%%%%%%%%%%%%%%%%%%%%%%%%%%
\prova One has $\e=\s\h^{\ppp}$, with $\s={\rm sign}(\e)$,
so that ${\hknu \F k \n}=
{\hknu \G k \n}=0$ for all $k<\ppp$ and all
$\n\in\ZZZ$, independently of the values of the constants
${\hknu \b 1 0},\ldots{\hknu \b {\ppp-1} 0}$. Moreover
${\hknu \bt k \n}={\hknu B k \n}=0$ for all $k<\ppp$.\qed
%%%%%%%%%%%%%%%%%%%%%%%%%%%%%%%%%%%%%%%%%%%%%%%%%%%%%%%%%%%%%%%%%%%%%%%%%

%%%%%%%%%%%%%%%%%%%%%%%%%%%%%%%%%%%%%%%%%%%%%%%%%%%%%%%%%%%%%%%%%%%%%%%%%
\begin{lemma}\label{lem:13}
The equations (\ref{eq:4.2}) are satisfied up to order $k=\ppp$,
for any choice of the constants
${\hknu \b 1 0},\ldots{\hknu \b {\ppp} 0}$.
\end{lemma}
%%%%%%%%%%%%%%%%%%%%%%%%%%%%%%%%%%%%%%%%%%%%%%%%%%%%%%%%%%%%%%%%%%%%%%%%%

%%%%%%%%%%%%%%%%%%%%%%%%%%%%%%%%%%%%%%%%%%%%%%%%%%%%%%%%%%%%%%%%%%%%%%%%%
\prova One has ${\hknu \G {\ppp} \,}=G(\a_{0}(t),A_{0},t+t_{0})$ and
${\hknu \F {\ppp} \,}=F(\a_{0}(t),A_{0},t+t_{0})$,
so that
\be\label{eq:4.4}
{\hknu \G {\ppp} \n}=\sum_{p\s_{0}+q\s'_{0}=\n}G_{\s_{0},\s'_{0}}(A_{0}),\qquad
{\hknu \F {\ppp} \n}=\sum_{p\s_{0}+q\s'_{0}=\n}F_{\s_{0},\s'_{0}}(A_{0}).
\ee
Thus, ${\hknu \bt {\ppp} \n}$ and ${\hknu B {\ppp} \n}$ can be
obtained from (\ref{eq:4.2}).
Finally ${\hknu \G {\ppp} 0}=M(t_{0})$ by definition, and one has
$M(t_{0})=0$ by Hypothesis \ref{hyp2}.
Hence also the last equation of (\ref{eq:4.2}) is satisfied.\qed
%%%%%%%%%%%%%%%%%%%%%%%%%%%%%%%%%%%%%%%%%%%%%%%%%%%%%%%%%%%%%%%%%%%%%%%%%
\salto
%%%%%%%%%%%%%%%%%%%%%%%%%%%%%%%%%%%%%%%%%%%%%%%%%%%%%%%%%%%%%%%%%%%%%%%%%

Let us set
\be\label{eq:4.5}\ba
& \hhh_{0} = \hhh^{(0)}\ppp^{(1)}\cdot\ldots\cdot\ppp^{(i_{0})},
&\quad&
 \sss_{0} = \sss^{(0)}\ppp^{(1)}\cdot\ldots\cdot\ppp^{(i_{0})},\\
& \hhh_{1} = \hhh_{0}+\hhh^{(1)}\ppp^{(2)}\cdot\ldots\cdot\ppp^{(i_{0})},
&\quad&
 \sss_{1} = \sss_{0}+\sss^{(1)}\ppp^{(2)}\cdot\ldots\cdot\ppp^{(i_{0})},\\
& \hhh_{2} = \hhh_{1}+\hhh^{(2)}\ppp^{(3)}\cdot\ldots\cdot\ppp^{(i_{0})},
&\quad&
 \sss_{2} = \sss_{1}+\sss^{(2)}\ppp^{(3)}\cdot\ldots\cdot\ppp^{(i_{0})}, \\
& \vdots
&\quad&
\vdots  \\
& \hhh_{i_{0}} = \hhh_{i_{0}-1}+\hhh^{(i_{0})},
&\quad&
 \sss_{i_{0}} = \sss_{i_{0}-1}+\sss^{(i_{0})}.
\ea\ee

%%%%%%%%%%%%%%%%%%%%%%%%%%%%%%%%%%%%%%%%%%%%%%%%%%%%%%%%%%%%%%%%%%%%%%%%%
\begin{lemma}\label{lem:14}
The equations (\ref{eq:4.2}) are satisfied up to
order $k=\ppp+\sss_{i_0}$ provided
${\hknu \b {\hhh_i} 0}=c_i$, with
$c_i$ the real root of a polynomial $P^{(i)}(c)$
of the $i$-th step of iteration step of the Newton-Puiseux process,
for $i=0,\ldots,i_0$, and ${\hknu \b {k'} 0}=0$
for all $k'\leq \hhh_{i_0}$, $k'\neq\hhh_i$ for any $i$.
\end{lemma}
%%%%%%%%%%%%%%%%%%%%%%%%%%%%%%%%%%%%%%%%%%%%%%%%%%%%%%%%%%%%%%%%%%%%%%%%%

%%%%%%%%%%%%%%%%%%%%%%%%%%%%%%%%%%%%%%%%%%%%%%%%%%%%%%%%%%%%%%%%%%%%%%%%%
\prova
If ${\hknu \b {k'} 0}=0$ for all $1<{k'}<\hhh_{0}$,
one has ${\hknu \G {k} 0}=0$ for all
$\ppp<k<\ppp+\sss_{0}$, while ${\hknu \G {\ppp+\sss_{0}} 0}=
P^{(0)}({\hknu \b {\hhh_0} 0})$,
so that ${\hknu \G {\ppp+\sss_{0}} 0}=0$ for
${\hknu \b {\hhh_0} 0}=c_0$.
Thus ${\hknu \G {k} 0}=0$ for
$\ppp+\sss_0<k<\ppp+\sss_{1}$ provided
${\hknu \b {k'} 0}=0$ for all $\hhh_0<{k'}<\hhh_{1}$,
while
${\hknu \G {\ppp+\sss_{1}} 0}=
P^{(1)}({\hknu \b {\hhh_1} 0})$,
so that ${\hknu \G {\ppp+\sss_{1}} 0}=0$ for
${\hknu \b {\hhh_1} 0}=c_1$, and so on.

Hence if we set
${\hknu \b {\hhh_i} 0}=c_i$, for all $i=0,\ldots,i_0$, and
${\hknu \b {k'} 0}=0$ for all ${k'}<\hhh_{i_0}$, $k'\neq\hhh_i$
for any $i=0,\ldots,i_0$, one has ${\hknu \G {k} 0}=0$ for all
$\ppp< k\le \ppp+\sss_{i_0}$.
Moreover, ${\hknu \F {k} \n}$ and ${\hknu \G {k} \n}$ are well-defined
for such values of $k$.
Hence (\ref{eq:4.2}) can be solved up to order $k=\ppp+\sss_{i_0}$,
indipendently of the values of the constants
${\hknu \b {k'} 0}$ for $k'>\hhh_{i_0}$.\qed
%%%%%%%%%%%%%%%%%%%%%%%%%%%%%%%%%%%%%%%%%%%%%%%%%%%%%%%%%%%%%%%%%%%%%%%%%
\salto
%%%%%%%%%%%%%%%%%%%%%%%%%%%%%%%%%%%%%%%%%%%%%%%%%%%%%%%%%%%%%%%%%%%%%%%%%

By Lemma \ref{lem:14} we can write
\be
\b_{0} = \b_{0}(\h) := c_0\h^{\hhh_0}+
c_1\h^{\hhh_1}+\ldots+
c_{i_0}\h^{\hhh_{i_0}}+\sum_{k\geq1}
{\hknu \b {\hhh_{i_0}+k}0}\h^{\hhh_{i_0}+k}.
\label{eq:4.6}\ee

%%%%%%%%%%%%%%%%%%%%%%%%%%%%%%%%%%%%%%%%%%%%%%%%%%%%%%%%%%%%%%%%%%%%%%%%%
\begin{lemma}\label{lem:15}
The equations (\ref{eq:4.2}) are satisfied up to any order
$k=\ppp+\sss_{i_0}+\k$, $\k\geq1$ provided the constants
${\hknu \b {\hhh_{i_0}+\k'} 0}$ are suitably fixed up to order $\k'=\k$.
\end{lemma}
%%%%%%%%%%%%%%%%%%%%%%%%%%%%%%%%%%%%%%%%%%%%%%%%%%%%%%%%%%%%%%%%%%%%%%%%%

%%%%%%%%%%%%%%%%%%%%%%%%%%%%%%%%%%%%%%%%%%%%%%%%%%%%%%%%%%%%%%%%%%%%%%%%%
\prova By substituting (\ref{eq:4.6}) and $\e=\s\h^{\ppp}$, with
$\s={\rm sign}(\e)$, in $\G_0(\e,\b_0)$ we obtain
\be\ba
\G_0(\s\h^{\ppp},\b_0(\h)) & = \s\h^{\ppp}\sum_{s_1,j\geq0}Q_{s_1,j}
\h^{s_1\ppp} \!\!\!\!\!\!\!
\sum_{\substack{m_0+\ldots+m_{i_0}+m=j\\ m,m_i\geq0}}
\!\!\!\!\!\!\!\!
J(j,m_{0},\ldots,m_{i_{0}},m) \, \times\\
& \times \,
\h^{m_0\hhh_0+\ldots+m_{i_0}\hhh_{i_0}}c_0^{m_0} \cdot
\ldots\cdot c_{i_0}^{m_{i_0}} \sum_{n\geq0}\h^{m\hhh_{i_0}+n}
\sum_{\substack{n_1+\ldots+n_{m}=n\\ n_i\geq1}}
{\hknu \b {\hhh_{i_0}+n_1}0}\ldots{\hknu \b {\hhh_{i_0}+n_{m}}0}
\label{eq:4.7}\ea\ee
where $Q_{s_1,j}=\FF_{s_1,j}^{(0)}\s^{s_1}$ and
\be
J(j,m_{0},\ldots,m_{i_{0}},m) := \frac{j!}{m_0!\ldots m_{i_0}!m!}.
\label{eq:4.9}\ee
For any $\k\geq1$ one has, by rearranging the sums,
\be
\s{\hknu \G {\ppp+\sss_{i_0}+\k} 0}=
\!\!\!\!\!\!\!\!\!\!\!\!\!\!
\!\!\!\!\!\!\!\!\!\!\!\!\!\!
\sum_{\substack{m,m_i,n,s_1,j\geq0\\ m_0+\ldots+m_{i_0}+m=j\\ 
s_1\ppp+m_0\hhh_0+\ldots+m_{i_0}\hhh_{i_0}+m\hhh_{i_0}=\sss_{i_0}+n}}
\!\!\!\!\!\!\!\!\!\!\!\!\!\!
\!\!\!\!\!\!\!\!\!\!\!\!\!\!
J(j,m_{0},\ldots,m_{i_{0}},m)\,
Q_{s_1,j}c_0^{m_0}\ldots c_{i_0}^{m_{i_0}}
\!\!\!\!\!\!\!\!\!\!\!\!\!
\sum_{\substack{n_1+\ldots+n_{m}=\k-n\\ n_{i}\geq1}}
\!\!\!\!\!\!\!\!\!\!\!\!\!
{\hknu \b {\hhh_{i_0}+n_1}0}\ldots{\hknu \b {\hhh_{i_0}+n_{m}}0},
\label{eq:4.8} \ee
so that the last equation of (\ref{eq:4.2}) gives for $\k\geq1$
\be \ba
& \sum_{\substack{m_i,s_1,j\geq0\\ m_0+\ldots+m_{i_0}+1=j \\
        s_1\ppp+m_0\hhh_0+\ldots+m_{i_0}\hhh_{i_0}=\sss_{i_0}}}
\!\!\!\!\!\!\!\!\!\!
\!\!\!\!\!\!\!\!\!\!
	m_{i_0}J(j,m_{0},\ldots,m_{i_{0}},m)\,
Q_{s_1,j}c_0^{m_0}\ldots c_{i_0-1}^{m_{i_0-1}}
c_{i_0}^{m_{i_0}-1}
        {\hknu \b {\hhh_{i_0}+\k}0} \\
& \qquad\qquad +
\!\!\!\!\!\!\!\!\!\!\!\!\!
\!\!\!\!\!\!\!\!\!\!\!\!\!
\sum_{\substack{m_i,s_1,j\geq0,m\geq2\\ m_0+\ldots+m_{i_0}+m=j \\
        s_1\ppp+m_0\hhh_0+\ldots+m_{i_0}
\hhh_{i_0}+m\hhh_{i_{0}}=\sss_{i_0}}}
\!\!\!\!\!\!\!\!\!\!\!\!\!
\!\!\!\!\!\!\!\!\!\!\!\!\!
        J(j,m_{0},\ldots,m_{i_{0}},m)\,
Q_{s_1,j}c_0^{m_0}\ldots c_{i_0}^{m_{i_0}}
\!\!\!\!\!\!\!\!\!
      \sum_{\substack{n_1+\ldots+n_{m}=\k\\ 1\leq n_{i}\leq\k-1}}
\!\!\!\!\!\!\!\!\!
      {\hknu \b {\hhh_{i_0}+n_1}0}\ldots{\hknu \b {\hhh_{i_0}+n_{m}}0} \\
& \qquad\qquad +
\!\!\!\!\!\!\!\!\!\!\!\!\!\!\!
\!\!\!\!\!\!\!\!\!\!\!\!\!\!
\sum_{\substack{m,m_i,s_1,j\geq0,n\geq1\\ m_0+\ldots+m_{i_0}+m=j\\ 
        s_1\ppp+m_0\hhh_0+\ldots+m_{i_0}\hhh_{i_0}+m\hhh_{i_{0}}=
\sss_{i_0}+n}}
\!\!\!\!\!\!\!\!\!\!\!\!\!\!\!\!
\!\!\!\!\!\!\!\!\!\!\!\!\!\!\!
	J(j,m_{0},\ldots,m_{i_{0}},m)\,
Q_{s_1,j}c_0^{m_0}\ldots c_{i_0}^{m_{i_0}}
\!\!\!\!\!\!\!\!\!\!\!
	\sum_{\substack{n_1+\ldots+n_{m}=\k-n\\ n_{i}\geq1}}
\!\!\!\!\!\!\!\!\!\!\!
	{\hknu \b {\hhh_{i_0}+n_1}0}\ldots{\hknu \b {\hhh_{i_0}+n_{m}}0}=0,
\ea
\label{eq:4.10}\ee
where all terms but those in the first line contain only coefficients
${\hknu \b {\hhh_{i_0}+\k'} 0}$ with $\k'<\k$.

Recall that by Hypothesis \ref{hyp3}
\be
\sum_{\substack{s_1,j\geq0\\ m_0+\ldots+m_{i_0}=j \\
        s_1\ppp+m_0\hhh_0+\ldots+m_{i_0}\hhh_{i_0}=\sss_{i_0}}}
\!\!\!\!\!\!\!\!\!\!
\!\!\!\!\!\!\!\!\!\!
	m_{i_0}J(j,m_{0},\ldots,m_{i_{0}},m)\,
Q_{s_1,j}c_0^{m_0}\ldots c_{i_0-1}^{m_{i_0-1}}
c_{i_0}^{m_{i_0}-1}
	= \frac{\de P^{(i_0)}}{\de c}
\left(c_{i_{0}}\right)=:C\neq0,
\label{eq:4.11} \ee
so that we can use (\ref{eq:4.10}) to express
${\hknu \b {\hhh_{i_0}+\k}0}$ in terms of the
coefficients ${\hknu \b {\hhh_{i_0}+\k'}0}$ of lower orders $\k'<\k$.
Thus we can conclude that the equations (\ref{eq:4.2})
are satisfied up to order ${k}$ provided the 
coefficients ${\hknu \b {\hhh_{i_0}+\k'}0}$ are fixed as
\be
{\hknu \b {\hhh_{i_0}+\k'}0}=-\frac{1}{C}
\tilde{G}^{[\k']}(c_0,\ldots,c_{i_0},{\hknu \b {\hhh_{i_0}+1}0},\ldots,
{\hknu \b {\hhh_{i_0}+\k'-1}0}),
\label{eq:4.12} \ee
for all $1\leq\k'\leq\k$, where $\tilde{G}^{[\k]}(c_0,\ldots,
c_{i_0},{\hknu \b {\hhh_{i_0}+1}0},\ldots,
{\hknu \b {\hhh_{i_0}+\k-1}0})$ is given by the sum of
the second and third lines in (\ref{eq:4.10}).\qed
%%%%%%%%%%%%%%%%%%%%%%%%%%%%%%%%%%%%%%%%%%%%%%%%%%%%%%%%%%%%%%%%%%%%%%%%%
\salto
%%%%%%%%%%%%%%%%%%%%%%%%%%%%%%%%%%%%%%%%%%%%%%%%%%%%%%%%%%%%%%%%%%%%%%%%%

We can summarise the results above into the following statement.

%%%%%%%%%%%%%%%%%%%%%%%%%%%%%%%%%%%%%%%%%%%%%%%%%%%%%%%%%%%%%%%%%%%%%%%%%
\begin{prop}\label{prop:2}
The equations (\ref{eq:4.2}) are satisfied to any order $k$ provided
the constants ${\hknu \b k 0}$ are suitably fixed.
In particular ${\hknu \bt k \n}={\hknu B k \n}={\hknu B k 0}=0$
for $k<\ppp$ and ${\hknu \b k 0}=0$ for $k<\hhh_{i_{0}}$, $k\neq\hhh_{i}$
for any $i=0,\ldots,i_{0}$.
\end{prop}
%%%%%%%%%%%%%%%%%%%%%%%%%%%%%%%%%%%%%%%%%%%%%%%%%%%%%%%%%%%%%%%%%%%%%%%%%

%%%%%%%%%%%%%%%%%%%%%%%%%%%%%%%%%%%%%%%%%%%%%%%%%%%%%%%%%%%%%%%%%%%%%%%%%
%%%%%%%%%%%%%%%%%%%%%%%%%%%%%%%%%%%%%%%%%%%%%%%%%%%%%%%%%%%%%%%%%%%%%%%%%
\zerarcounters
\section{Diagrammatic rules for the Puiseux series}
\label{sec:6}
%%%%%%%%%%%%%%%%%%%%%%%%%%%%%%%%%%%%%%%%%%%%%%%%%%%%%%%%%%%%%%%%%%%%%%%%%
%%%%%%%%%%%%%%%%%%%%%%%%%%%%%%%%%%%%%%%%%%%%%%%%%%%%%%%%%%%%%%%%%%%%%%%%%

\noindent In order to give a graphical representation of the coefficients
${\hknu \b k 0}$, ${\hknu \bt k \n}$ and ${\hknu B k \n}$
in (\ref{eq:4.1}), we shall consider a different tree expansion
with respect to that of Section \ref{sec:5}.
We shall perform an iterative construction, similar to the
one performed through the proof of Lemma \ref{lem:11},
starting from equations (\ref{eq:4.2}) for the coefficients
$\tilde{\b}_{\n}^{[k]}$, $B_{\n}^{[k]}$ for $k\geq\ppp$,
and from (\ref{eq:4.12}) for $\b_{0}^{[k]}$, $k\geq\hhh_{i_{0}}+1$.

Let us consider a tree with leaves.
We associate with each leaf $\eee$ a \emph{leaf label}
$\aaa_\eee=0,\ldots,i_{0}$.

For $k=\ppp$ we represent the coefficients
${\hknu {\tilde{\b}} \ppp \n}$ and ${\hknu B \ppp \n}$
as a line exiting a node, while for $k=\hhh_{i}$, $i=0,\ldots,i_{0}$
we represent ${\hknu \b {\hhh_{i}} 0}$ as a
line exiting a leaf with leaf label $\aaa_{i}$.

\begin{figure}[!ht]
\centering
\begin{minipage}[b]{7cm}
\centering
{
\psfrag{b}{${\hknu \b k 0}$}
\psfrag{btp}{${\hknu \bt k \n}$}
\psfrag{Bk}{${\hknu B k \n}$}
\psfrag{B}{$B$}
\psfrag{bt}{$\bt$}
\psfrag{b0}{$\b_{0}$}\psfrag{bt}{$\bt$}
\psfrag{[k]}{$[k]$}
\psfrag{+}{$+$}
\psfrag{2}{$2$}
\psfrag{1}{$1$}
\psfrag{0}{$0$}
\psfrag{n}{$\n$}
\includegraphics[width=11cm]{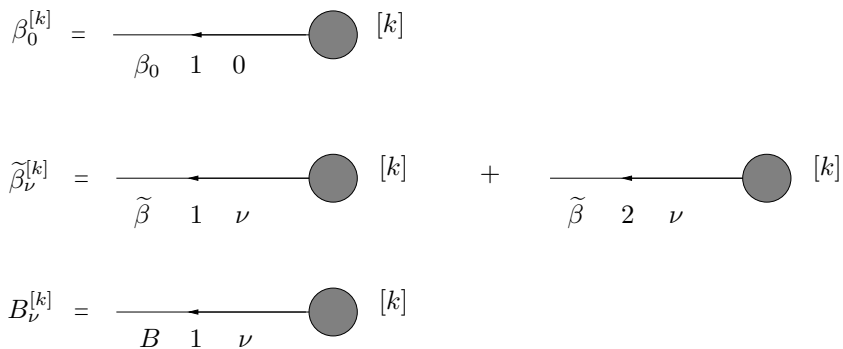}
}
\caption{\label{fig:4}\footnotesize{Graph elements.}}
\end{minipage}
\end{figure}

Now we represent each coefficient as a graph element according
to Figure \ref{fig:4}, as a line exiting a ball with
\emph{order label} $k$, with $k\geq\hhh_{i_{0}}+1$ for the
coefficients ${\hknu \b k 0}$, and
$k\geq\ppp+1$ for the coefficients ${\hknu {\tilde{\b}} k \n}$ and
${\hknu B k \n}$; we associate with the line a component
label $h_{\ell}\in\{\b_{0},\tilde{\b},B\}$, a degree
label $\d_{\ell}\in\{1,2\}$
with the constraint that $\d_{\ell}=1$ for $h_{\ell}=B,\b_0$,
and momentum label $\n_\ell\in\ZZZ$,
with the constraint that $\n_\ell\neq0$
for $h_\ell=\tilde{\b}$, while $\n_\ell=0$ for $h_\ell=\b_{0}$.

%%%%%%%%%%%%%%%%%%%%%%%%%%%%%%%%%%%%%%%%%%%%%%%%%%%%%%%%%%%%%%%%%%%%%%%%%
% figure 4
%%%%%%%%%%%%%%%%%%%%%%%%%%%%%%%%%%%%%%%%%%%%%%%%%%%%%%%%%%%%%%%%%%%%%%%%%
%%%%%%%%%%%%%%%%%%%%%%%%%%%%%%%%%%%%%%%%%%%%%%%%%%%%%%%%%%%%%%%%%%%%%%%%%

Hence we can represent the first three equations in (\ref{eq:4.2})
graphically, representing each factor $\b_{\n_{i}}^{[k_{i}]}$
and $B_{\n_{i}}^{[k_{i}]}$ in (\ref{eq:4.3}) as graph elements:
again the lines of such graph
elements enter the same node $\vvv_{0}$.

We associate with $\vvv_{0}$ a badge label
$b_{\vvv_{0}}\in\{0,1\}$ by setting
$b_{\vvv_{0}}=1$ for $h_{\ell_{0}}=\tilde{\b}$ and $\d_{\ell_{0}}=2$,
and for $h_{\ell_{0}}=B$ and $\n_{\ell_{0}}\neq0$.
We call $r_{\vvv_{0}}$ the number of the lines
entering $\vvv_{0}$ with component label $h=\b_{0},\tilde{\b}$, and
$s_{\vvv_{0}}$ the number of
the lines entering $\vvv_{0}$ with component label $h=B$, with the
constraint that if $b_{\vvv_{0}}=0$
one has $r_{\vvv_{0}}=0$ and $s_{\vvv_{0}}\geq2$. Finally we associate
with $\vvv_{0}$ two mode labels $\s_{\vvv_{0}},\s_{\vvv'_{0}}\in\ZZZ$
and the global mode label $\n_{\vvv_{0}}$ defined as
in (\ref{eq:5.1}), and we impose the conservation law
\be\label{eq:6.1}
\n_{\ell_{\vvv_{0}}}=\n_{\vvv_{0}}+\sum_{i=1}^{r_{\vvv_{0}}+
s_{\vvv_{0}}}\n_{\ell_{i}},
\ee
where $\ell_{1},\ldots,\ell_{r_{\vvv_{0}}+s_{\vvv_{0}}}$ are the
lines entering $\vvv_{0}$.

We also force the following conditions on the order labels
\be\label{eq:6.2}\ba
\sum_{i=1}^{r_{\vvv_{0}}+s_{\vvv_{0}}}k_{i}=k-\ppp,\qquad b_{\vvv_{0}}=1,\\
\sum_{i=1}^{s_{\vvv_{0}}}k_{i}=k,\qquad b_{\vvv_{0}}=0,
\ea\ee
which reflect the condition on the sums in (\ref{eq:4.3}).

Finally we associate with $\vvv=\vvv_{0}$ a
node factor $\NN^*_{\vvv}=\s^{b_{\vvv}}\NN_{\vvv}$, with $\s=\sign(\e)$
and $\NN_{\vvv}$ defined as in (\ref{eq:5.4}),
and with the line $\ell=\ell_{\vvv_0}$ a propagator $g^*_{\ell}=g_\ell$,
with $g_{\ell}$ defined as in (\ref{eq:5.3}).
The only difference with respect to Section \ref{sec:5} is
that the component label can assume the values $\tilde{\b},B$,
which have the r\^ole of $\a,A$ respectively.

The coefficients ${\hknu \b {\hhh_{i_0}+\k} 0}$, $\k\ge 1$,
have to be treated in a different way.

First of all we
point out that also the coefficients $\ol{\G}_0^{(k,j)}$ in (\ref{eq:2.20})
can be represented in
terms of sum of trees with leaves as in Section \ref{sec:5}.
In fact we can repeat the iterative construction
of Lemma \ref{lem:11}, simply by defining $\TT_{k,0,\G,j}$ as
the set of the trees contributing to $\ol{\G}_0^{(k,j)}$,
setting $g_{\ell_0}=1$,
$h_{\ell_0}=\G$, $\d_{\ell_0}=1$, $\n_{\ell_0}=0$, $b_{\vvv_0}=1$ and
\be\label{eq:6.3}
\NN_{\vvv_0}=
\frac{(i\s_{\vvv_0})^{r_{\vvv_0}}\dpr_A^{s_{\vvv_0}}}{r_{\vvv_0}!s_{\vvv_0}!}
G_{\s_{\vvv_0},\s_{\vvv_0}'}(A_0,t_0),
\ee
and no further difficulties arise.

\begin{figure}[!ht]
\centering
\begin{minipage}[b]{10cm}
\centering{
\psfrag{b0}{$\b_{0}$}
\psfrag{s1+1}{$s_{1}+1$}
\psfrag{s0}{$s_{0}$}
\psfrag{[k1]}{$[k_{1}]$}
\psfrag{[ks0]}{$[k_{s_{0}'}]$}
\psfrag{0}{$0$}
\psfrag{1}{$1$}
\psfrag{nodes}{nodes}
\psfrag{leaves}{leaves}
\includegraphics[width=8cm]{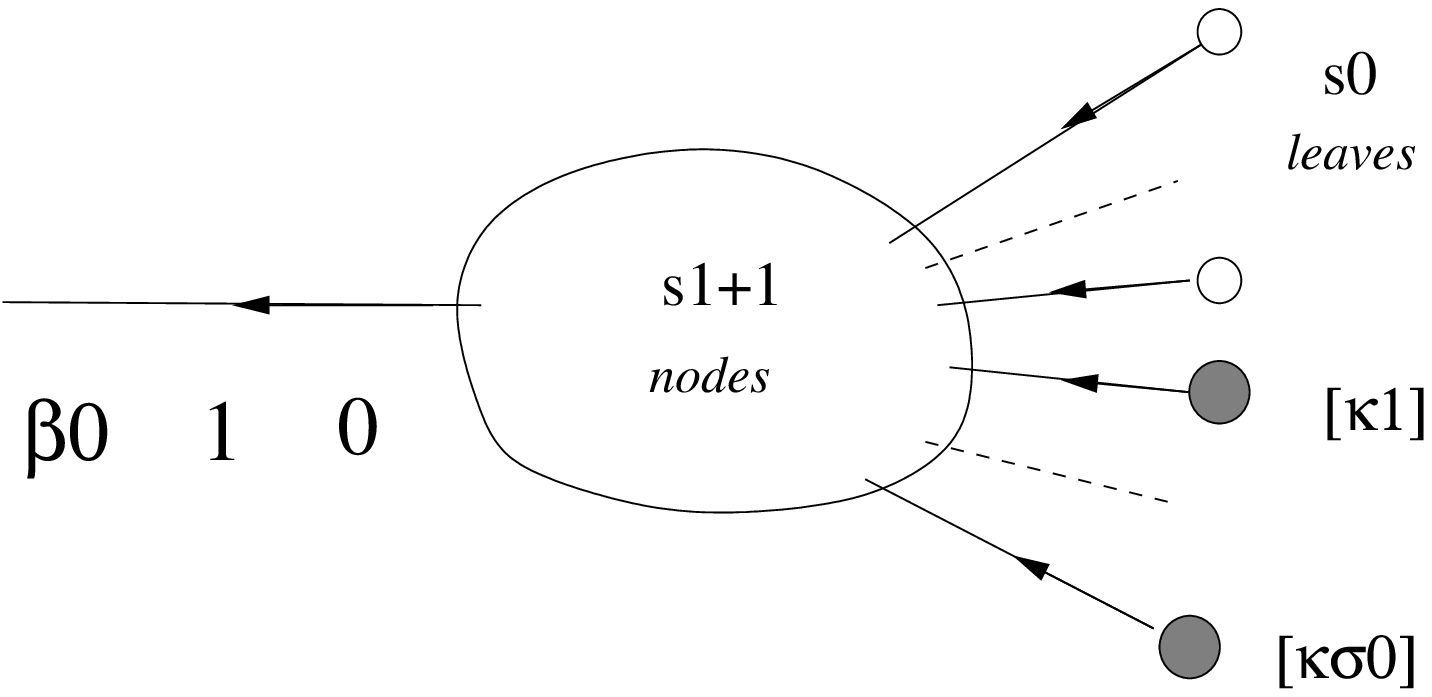}
}
\caption{\label{fig:5}
\footnotesize{A tree contributing to ${\hknu \b k 0}$.}}
\end{minipage}
\end{figure}

Recall that the
coefficients $Q_{s_{1},j}$ in (\ref{eq:4.10}) are defined as
$Q_{s_{1},j}=\FF_{s_{1},j}^{(0)}\s^{s_{1}}=
{\eknu{\ol{\G}}{s_{1}+1,j}0}\s^{s_{1}}$
so that
\be\label{eq:6.4}
Q_{s_{1},j}=\s^{s_{1}}
\!\!\!
\sum_{\th\in\TT_{s_{1}+1,0,\G,j}}
\!\!\!
\Val(\th).
\ee
Hence the summands in 
the second and third lines in
(\ref{eq:4.10}) can be
imagined as ``some'' of the trees in $\TT_{s_{1}+1,0,\G,j}$
where ``some'' leaves are substituted by graph elements
with $h_\ell=\b_{0}$. More precisely we shall consider
only trees $\th$ of the form depicted in Figure \ref{fig:5},
with $s_{1}+1$ nodes, $s_{0}=s_{0,0}+\ldots+s_{0,i_0}$
leaves, where $s_{0,\aaa}$ is the number of the leaves with
leaf label $\aaa$,
and $s_{0}'$ graph elements with
$h_\ell=\b_{0}$, such that
\be\label{eq:6.5}\ba
&s_{1}\ppp+\sum_{i=0}^{i_0}s_{0,i}\hhh_i+s_{0}'\hhh_{i_{0}}=
            \sss_{i_{0}}+n,\\
&\sum_{i=1}^{s_{0}'}k_{i}=(s_{0}'-1)\hhh_{i_{0}}+k-n,
\ea\ee
for a suitable $0\leq n\leq k-\hhh_{i_{0}}$, with the constraint that
when $n=0$ one has $s_{0}'\geq2$. We shall call $\ell_{i}$ the $s_{0}'$
lines with $h_{\ell_{i}}=\b_{0}$. Such conditions
express the condition on the sums
in the second and third lines
in (\ref{eq:4.10}).

%%%%%%%%%%%%%%%%%%%%%%%%%%%%%%%%%%%%%%%%%%%%%%%%%%%%%%%%%%%%%%%%%%%%%%%%%
% figure 5
%%%%%%%%%%%%%%%%%%%%%%%%%%%%%%%%%%%%%%%%%%%%%%%%%%%%%%%%%%%%%%%%%%%%%%%%%
%%%%%%%%%%%%%%%%%%%%%%%%%%%%%%%%%%%%%%%%%%%%%%%%%%%%%%%%%%%%%%%%%%%%%%%%%

The propagators of the lines exiting any among the $s_{1}+1$ nodes
and the node factors of the nodes (except the root line
and the node which the root line exits) are
$g^*_{\ell}=g_\ell$ and $\NN^*_\vvv=\s^{b_{\vvv}}\NN_\vvv$
with the component labels assuming the values $\tilde{\b},B$,
which have the r\^ole of $\a,A$, respectively.
We associate with the root line a propagator
\be\label{eq:6.6}
g^*_{\ell_{0}}=- \frac{1}{C},
\ee
where $C$ is defined in (\ref{eq:4.11}), while the node $\vvv=\vvv_0$
which the root line exits will have a node factor
$\NN_{\vvv_0}^*=\NN_{\vvv_0}$ as in (\ref{eq:6.3})
Finally we associate with each leaf $\eee$ a \emph{leaf factor}
${\NN}^*_\eee=c_{\aaa_\eee}$.

We now iterate such a process until only nodes or leaves appear.
We shall call \emph{allowed trees} all the trees obtained
in such a recursive way, and
we shall denote with $\Th_{k,\n,h}$ the set of allowed trees
with order $k$, total momentum $\n$
and total component label $h$.

Given an allowed tree $\th$ we denote
with $N(\th)$, $L(\th)$ and $E(\th)$ the set of nodes,
lines and leaves of $\th$ respectively,
and we denote with ${E}_{\aaa}(\th)$ the set of leaves in $\th$
with leaf label $\aaa$.
We point out that $E(\th)=E_{0}(\th)\amalg\ldots\amalg E_{i_{0}}(\th)$.
We shall define the \emph{value} of $\th$ as
\be\label{eq:6.7}
\Val^*(\th)=\left(\prod_{\ell\in L(\th)}g_\ell^*\right)
\left(\prod_{\vvv\in N(\th)}\NN_\vvv^*\right)
\left(\prod_{\eee\in E(\th)}\NN_\eee^*\right).
\ee

Finally, we denote with $\L(\th)$ the set of the lines
(exiting a node) in $\th$ with component label $h=\b_{0}$ and with
$N^*(\th)$ the nodes with $b_\vvv=1$; then we
associate with each node in $N^*(\th)$, with each leaf
in ${E}_\aaa(\th)$ and with each line in $\L(\th)$ a
\emph{weight} $\ppp$, $\hhh_\aaa$ and $\hhh_{i_{0}}-\ppp-\sss_{i_{0}}$,
respectively, and we call \emph{order} of $\th$ the number
\be\label{eq:6.8}
k(\th)=\ppp|N^*(\th)|+(\hhh_{i_{0}}-\ppp-\sss_{i_{0}})|\L(\th)|+
\sum_{\aaa=0}^{i_{0}}\hhh_\aaa|E_\aaa(\th)|.
\ee
Note that $\hhh_{i_{0}}-\ppp-\sss_{i_{0}}<0$.

%%%%%%%%%%%%%%%%%%%%%%%%%%%%%%%%%%%%%%%%%%%%%%%%%%%%%%%%%%%%%%%%%%%%%%%%%
\begin{lemma}\label{lem:16}
The Fourier coefficients ${\hknu {\b} k 0}$, ${\hknu \bt k \n}$
and ${\hknu {{B}} {k} \n}$ can be written in terms of trees as
\begin{subequations}
\begin{align}
{\hknu {\b} {k} 0} & = \sum_{\th\in\Th_{k,0,\b_{0}}}
\Val^*(\th),\qquad k\geq\hhh_{i_0}+1,
\label{eq:6.9a} \\
{\hknu \bt {k} \n} &= \sum_{\th\in\Th_{k,\n,\bt}}
\Val^*({\th}),\qquad k\geq\ppp,
\label{eq:6.9b} \\
{\hknu {B} k \n} & =\sum_{\th\in\Th_{k,\n,B}}
\Val^*(\th),\qquad k\geq\ppp.
\label{eq:6.9c}
\end{align}
\label{eq:6.9} \end{subequations}
%
%\vskip-.3truecm \noindent 
\end{lemma}
%%%%%%%%%%%%%%%%%%%%%%%%%%%%%%%%%%%%%%%%%%%%%%%%%%%%%%%%%%%%%%%%%%%%%%%%%

%%%%%%%%%%%%%%%%%%%%%%%%%%%%%%%%%%%%%%%%%%%%%%%%%%%%%%%%%%%%%%%%%%%%%%%%%
\prova We only have to prove that an allowed tree contributing
to the Fourier coefficients ${\hknu \b k 0}$,
${\hknu \bt k \n}$ and ${\hknu {{B}} {k} \n}$ has order $k$.
We shall perform the proof by induction on $k\geq\hhh_{i_{0}}+1$
for the coefficients ${\hknu \b k 0}$, and
$k\geq\ppp$ for ${\hknu \bt k \n}$ and ${\hknu B k \n}$.
Let us set $f_{\tilde{\b}}=\tilde{\b}$ and $f_{B}=B$.
An allowed tree $\th$ contributing to
${\hknu f \ppp {h,\n}}$ has
only one node so that $k(\th)=\ppp$, while an allowed
tree $\ol{\th}$ contributing to ${\hknu \b {\hhh_{i}+1} 0}$,
$i=0,\ldots,i_{0}$ has $s_1+1$ nodes, $s_0=s_{0,0}+\ldots+s_{0,i_0}$
leaves, and one line $\L(\ol{\th})$, and \emph{via} the conditions
(\ref{eq:6.5}) we have
$s_1\ppp+s_{0,0}\hhh_0+\ldots+s_{0,i_0}\hhh_{i_0}=\sss_{i_0}+1$,
so that $k(\ol{\th})=\hhh_{i_0}+1$.

Let us suppose first that for all $k'<k$, an allowed
tree $\th'$ contributing to
${\hknu \b {k'} 0}$ has order $k(\th')=k'$.
By the inductive hypothesis, the order of a tree $\th$ contributing
to ${\hknu \b k 0}$ is (we refer again to Figure \ref{fig:5} for notations)
\be\label{eq:6.10}
k(\th)=(s_{1}+1)\ppp+\sum_{i=0}^{i_0}s_{0,i}\hhh_{i}+
\sum_{i=1}^{s_{0}'}k_{i}+\hhh_{i_{0}}-\ppp-\sss_{i_{0}},
\ee
and \emph{via} the conditions in (\ref{eq:6.5}) we obtain $k(\th)=k$.

\begin{figure}[!ht]
\centering
\begin{minipage}[b]{10cm}
\centering{
\psfrag{v0}{$\vvv_{0}$}
\psfrag{s1}{$s_{1}$}
\psfrag{s0}{$s_{0}$}
\psfrag{s0'}{$s_{0}'$}
\psfrag{s1'}{$s_{1}'$}
\psfrag{q}{$\th$}
\includegraphics[width=7cm]{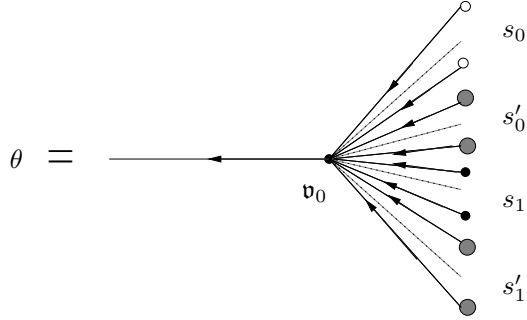}
}
\caption{\label{fig:6}
\footnotesize{An allowed tree contributing
to ${\hknu f \k {h,\n}}$.}}
\end{minipage}
\end{figure}

Let us suppose now that the inductive hypothesis holds
for all trees $\th'$ contributing to ${\hknu f {k'} {h,\n}}$, $k'<k$.
An allowed tree $\th$ contributing to ${\hknu f k {h,\n}}$
is of the form depicted in Figure \ref{fig:6}, where $s_{0,\aaa}$
is the number of the lines exiting a leaf
with leaf label $\aaa$ and entering $\vvv_{0}$, $s_{1}$ is the
number of the lines exiting a node and entering $\vvv_{0}$,
and $s_{0}',s_{1}'$ are the graph elements entering $\vvv_{0}$ with
component label $\b_{0}$ and either $\tilde{\b}$ or $B$, respectively.

%%%%%%%%%%%%%%%%%%%%%%%%%%%%%%%%%%%%%%%%%%%%%%%%%%%%%%%%%%%%%%%%%%%%%%%%%
% figure 6
%%%%%%%%%%%%%%%%%%%%%%%%%%%%%%%%%%%%%%%%%%%%%%%%%%%%%%%%%%%%%%%%%%%%%%%%%
%%%%%%%%%%%%%%%%%%%%%%%%%%%%%%%%%%%%%%%%%%%%%%%%%%%%%%%%%%%%%%%%%%%%%%%%%

If $b_{\vvv_{0}}=1$, by the inductive hypothesis the order of such
a tree is given by
\be\label{eq:6.11}
k(\th)=(s_{1}+1)\ppp+\sum_{i=0}^{i_0}s_{0,i}\hhh_{i}+
\sum_{i=1}^{s_{0}'+s_{1}'}k_{i},
\ee
and by the first condition in (\ref{eq:6.2}) we have
$k(\th)=k$. Otherwise if $b_{\vvv_{0}}=0$, we have $s_{0}+s_{0}'=0$
and, by the inductive hypothesis,
\be\label{eq:6.12}
k(\th)=s_{1}\ppp+\sum_{i=1}^{s_{1}'}k_{i}=k,
\ee
\emph{via} the second condition in (\ref{eq:6.2}).\qed
%%%%%%%%%%%%%%%%%%%%%%%%%%%%%%%%%%%%%%%%%%%%%%%%%%%%%%%%%%%%%%%%%%%%%%%%%

%%%%%%%%%%%%%%%%%%%%%%%%%%%%%%%%%%%%%%%%%%%%%%%%%%%%%%%%%%%%%%%%%%%%%%%%%
\begin{lemma}\label{lem:17}
Let $\qqq:=\min\{\hhh_{{0}},\ppp\}$ and let us define
\be\label{eq:6.13}
M=2 \frac{\sss_{i_{0}}}{\qqq}+3.
\ee
Then for all $\th\in\Th_{k,\n,h}$ one has
\be\label{eq:6.14}
|L(\th)|\leq Mk.
\ee
\end{lemma}
%%%%%%%%%%%%%%%%%%%%%%%%%%%%%%%%%%%%%%%%%%%%%%%%%%%%%%%%%%%%%%%%%%%%%%%%%

As the proof is rather technical we shall perform it
in Appendix \ref{app:B}.

The convergence of the series (\ref{eq:4.1}) for small $\h$
follows from the following result.

%%%%%%%%%%%%%%%%%%%%%%%%%%%%%%%%%%%%%%%%%%%%%%%%%%%%%%%%%%%%%%%%%%%%%%%%%
\begin{prop}\label{prop:3}
The formal solution (\ref{eq:4.1}) of the system (\ref{eq:2.9}),
given by the recursive
equations (\ref{eq:4.2}) and (\ref{eq:4.12}),
converges for $\h$ small enough.
\end{prop}
%%%%%%%%%%%%%%%%%%%%%%%%%%%%%%%%%%%%%%%%%%%%%%%%%%%%%%%%%%%%%%%%%%%%%%%%%

%%%%%%%%%%%%%%%%%%%%%%%%%%%%%%%%%%%%%%%%%%%%%%%%%%%%%%%%%%%%%%%%%%%%%%%%%
\prova By Lemma \ref{lem:17}, the number of unlabeled trees
of order $k$ is bounded by $4^{Mk}$. Thus, the sum over
all labels except the mode labels and the momenta is bounded
by a constant to the power $k$ because all such labels can
assume only a finite number of values. The bound
for each node factor is the same as in Proposition \ref{prop:1},
while the propagators
can be bounded by
\be\label{eq:6.15}
|g_{\ell}^*|\leq\max\left\{\left| \frac{\o'(A_{0})}{\o^2}\right|,
\left|\frac{1}{\o'(A_{0})}\right|,\left|\frac{1}{\o}\right|,
\left|\frac{1}{C}\right|,1\right\},
\ee
so that the product over all the lines can be bounded again 
by a constant to the power $k$.
The product over the leaves factors is again bounded 
by a constant to the power $k$,
while the sum over the mode labels which uniquely 
determine the momenta can be performed
by using for each node half the exponential decay factor 
provided by (\ref{eq:5.10}). Thus we obtain
\be\label{eq:6.16}
|\tilde{\b}^{[k]}_{\n}|\leq C_{1}C^{k}_{2}e^{-\k|\n|/2},\quad
|{B}^{[k]}_{\n}|\leq C_{1}C^{k}_{2}e^{-\k|\n|/2},
\ee
for suitable constants $C_{1}$ and $C_{2}$.
Hence we obtain the convergence for the series (\ref{eq:4.1}),
for $|\h|\leq C_{2}^{-1}$. \qed
%%%%%%%%%%%%%%%%%%%%%%%%%%%%%%%%%%%%%%%%%%%%%%%%%%%%%%%%%%%%%%%%%%%%%%%%%
\salto
%%%%%%%%%%%%%%%%%%%%%%%%%%%%%%%%%%%%%%%%%%%%%%%%%%%%%%%%%%%%%%%%%%%%%%%%%

The discussion above ends the proof of Theorem \ref{thm:2}. 

%%%%%%%%%%%%%%%%%%%%%%%%%%%%%%%%%%%%%%%%%%%%%%%%%%%%%%%%%%%%%%%%%%%%%%%%%
%%%%%%%%%%%%%%%%%%%%%%%%%%%%%%%%%%%%%%%%%%%%%%%%%%%%%%%%%%%%%%%%%%%%%%%%%
\zerarcounters
\section{Higher order subharmonic Melnikov functions}
\label{sec:7}
%%%%%%%%%%%%%%%%%%%%%%%%%%%%%%%%%%%%%%%%%%%%%%%%%%%%%%%%%%%%%%%%%%%%%%%%%
%%%%%%%%%%%%%%%%%%%%%%%%%%%%%%%%%%%%%%%%%%%%%%%%%%%%%%%%%%%%%%%%%%%%%%%%%

\noindent Now we shall see how to extend the results above
when the Melnikov function vanishes identically.

We are searching for a solution of the form $(\a(t),A(t))$ with
$\a(t)=\a_{0}(t)+\b_{0}+\bt(t)$ and $A(t)=A_{0}+B(t)$, where
\be \bt(t)=\sum_{\substack{\n\in\ZZZ\\ \n\neq0}}
e^{i\n\o t}\b_\n(\e,\b_{0}), \qquad
B(t)=\sum_{\n\in\ZZZ}e^{i\n\o t}B_\n(\e,\b_{0}).
\label{eq:7.1} \ee

First of all, we notice that we can formally write the equations of motion as
\be\label{eq:7.2}\left\{\ba
&\ol{\b}^{(k)}_{\n}(\b_{0})=
\frac{\ol{\F}^{(k)}_{\n}(\b_{0})}{i\o\n}+
\o'(A_{0}) \frac{\ol{\G}^{(k)}_{\n}(\b_{0})}{(i\o\n)^2}, \qquad
\ol{B}^{(k)}_{\n}(\b_{0})= \frac{\ol{\G}^{(k)}_{\n}(\b_{0})}{ i\o\n},
\qquad \n\neq0 \\
&\ol{B}^{(k)}_{0}(\b_{0})=-\frac{\ol{\F}^{(k)}_{0}(\b_{0})}{\o'(A_{0})}.
\ea\right.\ee
where the notations in (\ref{eq:2.22}) have been used,
up to any order $k$, provided
\be\label{eq:7.3}
\ol{\G}_{0}(\e,\b_{0})=0.
\ee

If $M(t_{0})$ vanishes identically, by (\ref{eq:2.24}) we
have $\ol{\G}_{0}^{(1,j)}=0$ for all $j\geq0$, that is
$\ol{\G}^{(1)}_{0}(\b_{0})=0$, for all $\b_{0}$, and
hence $\ol{\G}_{0}(\e,\b_{0})=\e^2\FF^{(2)}(\e,\b_{0})$, with
$\FF^{(2)}$ a suitable function analytic in $\e$, $\b_{0}$.

Thus, we can solve the equations of motion up to
the first order in $\e$, and the parameter
$\b_{0}$ is left undetermined. More precisely we obtain
\be
\b_\n=\e\ol{\b}^{(1)}_\n+\e\tilde{\b}_\n^{(1)}(\e,\b_{0}), \qquad
B_\n=\e \ol{B}^{(1)}_\n+\e\tilde{B}_\n^{(1)}(\e,\b_{0}), 
\label{eq:7.4} \ee
where $\ol{\b}_\n^{(1)}$, $\ol{B}_\n^{(1)}$ solve the
equation of motion up to the first order in $\e$,
while $\tilde{\b}_\n^{(1)}$, $\tilde{B}_\n^{(1)}$
are the corrections to be determined.

Now, let us set
\be\label{eq:7.5}
M_{0}(t_{0})=M(t_{0}),\qquad
M_{1}(t_{0})=\ol{\G}_{0}^{(2)}(0,t_{0}),
\ee
where $\ol{\G}^{(k)}_\n(\b_{0},t_{0})=\ol{\G}^{(k)}_\n(\b_{0})$,
\ie we are stressing
the dependence of ${\eknu {\ol{\G}} {k,j} \n}$ on $t_{0}$.
We refer to $M_{1}(t_{0})$
as the \emph{second order subharmonic Melnikov function}.
Notice that $M_{0}(t_{0})=\ol{\G}_{0}^{(1)}(0,t_{0})$.

If there exist $t_{0}\in[0,2\p)$ and $\nnn_{1}\in\NNN$ such that
$t_{0}$ is a zero of order $\nnn_{1}$ for the second order
subharmonic Melnikov function, that is
\be\label{eq:7.6}
\frac{\de^{k}}{\de t_{0}^{k}}M_{1}(t_{0})=0\;\;\;\forall\;0
\leq k\leq\nnn_{1}-1,\phantom{ and }
D=D(t_{0}):=\frac{\de^{\nnn_{1}}}{\de t_{0}^{\nnn_{1}}}M_{1}(t_{0})\neq0,
\ee
then we can repeat the analysis of the previous Sections 
to obtain the existence of a subharmonic
solution. In fact, we have
\be\label{eq:7.7}
\FF^{(2)}(\e,\b_{0}):=\sum_{k,j\geq0}\e^{k}\b_{0}^{j}\FF_{k,j}^{(2)},
\qquad \FF_{k,j}^{(2)}=\ol{\G}^{(k+2,j)}_{0}(t_{0}),
\ee
where $t_{0}$ has to be fixed as the zero of $M_{1}(t_{0})$,
so that, as
\be\label{eq:7.8}
(-\o(A_{0}))^{-j} \frac{\de^{j}}{\de t_{0}^{j}}M_{1}(t_{0})=
j!\ol{\G}_{0}^{(2,j)}(t_{0}),
\ee
for all $j$, as proved in \cite{GBD} with a different notation,
we can construct the Newton polygon of $\FF^{(2)}$,
which is $\b_{0}$-general of order $\nnn_{1}$ by (\ref{eq:7.6}),
to obtain $\tilde{\b}^{(1)}$, $\tilde{B}^{(1)}$ and $\b_{0}$
as Puiseux series in $\e$, provided at each step of the
iteration of the Newton-Puiseux algorithm one has a real root.

Otherwise, if $M_{1}(t_{0})$ vanishes identically, we have
$\ol{\G}_{0}^{(2)}(\b_{0})=0$ for all $\b_{0}$,
so that we can solve the equations of motion up to
the second order in $\e$ and the parameter
$\b_{0}$ is still undetermined. Hence we set
$M_{2}(t_{0})=\ol{\G}_{0}^{(3)}(0,t_{0})$ and so on.

In general if $M_{k'}(t_{0})\equiv0$, for all $k'=0,\ldots,\k-1$,
we have $\ol{\G}_{0}(\e,\b_{0})=\e^{k'}\FF^{(k')}(\e,\b_{0})$,
so that we can solve the equations of motion up to the $\k$-th
order in $\e$, and obtain
\begin{subequations}
\begin{align}
\b_\n & = \e\ol{\b}^{(1)}_\n+\ldots+\e^{\k}\ol{\b}^{(\k)}_\n+
\e^{\k}\tilde{\b}_\n^{(k)}(\e,\b_{0}),
\label{eq:7.9a} \\
B_\n & = \e \ol{B}^{(1)}_\n+\ldots+\e^{\k}\ol{B}^{(\k)}_\n+
\e^{\k}\tilde{B}_\n^{(k)}(\e,\b_{0}),
\label{eq:7.9b}
\end{align}
\label{eq:7.9} \end{subequations}
\vskip-.3truecm
\noindent where $\ol{\b}_\n^{(k')}$, $\ol{B}_\n^{(k')}$,
$k'=0,\ldots,\k-1$ solve the equation of motion up to the
$\k$-th order in $\e$, while $\tilde{\b}_\n^{(\k)}$,
$\tilde{B}_\n^{(\k)}$ are the correction to be determined.

Hence we can weaken Hypotheses \ref{hyp2} and \ref{hyp3} as follows.

%%%%%%%%%%%%%%%%%%%%%%%%%%%%%%%%%%%%%%%%%%%%%%%%%%%%%%%%%%%%%%%%%%%%%%%%%
\begin{hyp}\label{hyp4}
There exists $\k\geq0$ such that for all $k'=0,\ldots,\k-1$,
$M_{k'}(t_{0})$ vanishes identically, and there exist
$t_{0}\in[0,2\p)$ and $\nnn\in\NNN$ such that
\be\label{eq:7.11}
\frac{\de^{j}}{\de t_{0}^{j}}M_\k(t_{0})=0\;\;\;\forall\;0
\leq j\leq\nnn-1,\phantom{ and }
D=D(t_{0}) := \frac{\de^{\nnn}}{\de t_{0}^{\nnn}}M_\k(t_{0})\neq0 ,
\ee
that is $t_{0}$ is a zero of order $\nnn$ for the $\k$-th
order subharmonic Melnikov function.
\end{hyp}
%%%%%%%%%%%%%%%%%%%%%%%%%%%%%%%%%%%%%%%%%%%%%%%%%%%%%%%%%%%%%%%%%%%%%%%%%

%%%%%%%%%%%%%%%%%%%%%%%%%%%%%%%%%%%%%%%%%%%%%%%%%%%%%%%%%%%%%%%%%%%%%%%%%
\begin{hyp}\label{hyp5}
There exists $i_{0}\geq0$ such that at the $i_{0}$-th step of
the iteration of the Newton-Puiseux algorithm for $\FF^{(\k)}$,
there exists a polynomial $P^{(i_{0})}=P^{(i_{0})}(c)$ which has a
simple root $c^*\in\RRR$.
\end{hyp}
%%%%%%%%%%%%%%%%%%%%%%%%%%%%%%%%%%%%%%%%%%%%%%%%%%%%%%%%%%%%%%%%%%%%%%%%%

Thus we have the following result.

%%%%%%%%%%%%%%%%%%%%%%%%%%%%%%%%%%%%%%%%%%%%%%%%%%%%%%%%%%%%%%%%%%%%%%%%%
\begin{teo}\label{thm:3}
Consider a periodic solution with frequency $\o=p/q$ for
the system (\ref{eq:2.1}), and assume that Hypotheses \ref{hyp1},
\ref{hyp4} and \ref{hyp5} are satisfied.
Then there exists an
explicitly computable value $\e_{0}>0$ such that for $|\e|<\e_{0}$
the system (\ref{eq:2.1}) has at least
one subharmonic solution of order $q/p$.
Such a solution admits a convergent power series in $|\e|^{1/\nnn!}$,
and hence a convergent Puiseux series in $|\e|$.
\end{teo}
%%%%%%%%%%%%%%%%%%%%%%%%%%%%%%%%%%%%%%%%%%%%%%%%%%%%%%%%%%%%%%%%%%%%%%%%%

The proof can be easily obtained suitably modifying the proof
of Theorem \ref{thm:2}.

Now, call $\Re_n^{(\k)}$ the set of real roots of the polynomials
obtained at the $n$-th step of iteration of the Newton-Puiseux
process for $\FF^{(\k)}$. Again if $\nnn$ is even
we can not say \emph{a priori} whether a formal solution
exists at all. However, if $\Re_n^{(\k)}\neq\emptyset$
for all $n\geq0$, then we obtain a convergent
Puiseux series as in Section \ref{sec:3}.

Finally, as a corollary, we have the following result.

%%%%%%%%%%%%%%%%%%%%%%%%%%%%%%%%%%%%%%%%%%%%%%%%%%%%%%%%%%%%%%%%%%%%%%%%%
\begin{teo}\label{thm:4}
Consider a periodic solution with frequency $\o=p/q$
for the system (\ref{eq:2.1}). Assume that Hypotheses \ref{hyp1}
and \ref{hyp4} are satisfied with $\nnn$ odd. Then for $\e$ small
enough the system (\ref{eq:2.1})
has at least one subharmonic solution of order $q/p$.
Such a solution admits a convergent power series in $|\e|^{\nnn!}$,
and hence a convergent Puiseux series in $\e$.
\end{teo}
%%%%%%%%%%%%%%%%%%%%%%%%%%%%%%%%%%%%%%%%%%%%%%%%%%%%%%%%%%%%%%%%%%%%%%%%%

Again the proof is a suitable modification of the
proof of Theorem \ref{thm:1}.

%%%%%%%%%%%%%%%%%%%%%%%%%%%%%%%%%%%%%%%%%%%%%%%%%%%%%%%%%%%%%%%%%%%%%%%%%
\salto\salto
%%%%%%%%%%%%%%%%%%%%%%%%%%%%%%%%%%%%%%%%%%%%%%%%%%%%%%%%%%%%%%%%%%%%%%%%%
\noindent{\bf Acknowledgements.} We thank Edoardo Sernesi
for useful discussions.
%%%%%%%%%%%%%%%%%%%%%%%%%%%%%%%%%%%%%%%%%%%%%%%%%%%%%%%%%%%%%%%%%%%%%%%%%

%%%%%%%%%%%%%%%%%%%%%%%%%%%%%%%%%%%%%%%%%%%%%%%%%%%%%%%%%%%%%%%%%%%%%%%%%
\appendix
%%%%%%%%%%%%%%%%%%%%%%%%%%%%%%%%%%%%%%%%%%%%%%%%%%%%%%%%%%%%%%%%%%%%%%%%%

%%%%%%%%%%%%%%%%%%%%%%%%%%%%%%%%%%%%%%%%%%%%%%%%%%%%%%%%%%%%%%%%%%%%%%%%%
%%%%%%%%%%%%%%%%%%%%%%%%%%%%%%%%%%%%%%%%%%%%%%%%%%%%%%%%%%%%%%%%%%%%%%%%%
\zerarcounters
\section{On the genericity of Hypothesis 3}
\label{app:A}
%%%%%%%%%%%%%%%%%%%%%%%%%%%%%%%%%%%%%%%%%%%%%%%%%%%%%%%%%%%%%%%%%%%%%%%%%
%%%%%%%%%%%%%%%%%%%%%%%%%%%%%%%%%%%%%%%%%%%%%%%%%%%%%%%%%%%%%%%%%%%%%%%%%

\noindent Here we want to show that Hypothesis \ref{hyp3} is generic
on the space of the coefficients of the polynomials.
More precisely, we shall show that
given a polynomial of the form
\be\label{a.1}
P(a,c)=\sum_{i=0}^{n}a_{n-i}c^{i},\qquad n\geq1,\qquad
a:=(a_{0},\ldots,a_n),
\ee
the set of parameters $(a_{0},\ldots,a_n)\in\RRR^{n+1}$
for which $P(a,c)$ has multiples roots, is a proper
Zariski-closed\footnote{See for instance \cite{SHAF}.}
subset of $\RRR^{n+1}$. Notice that a polynomial $P=P(a,c)$
has a multiple root $c^*$ if and only if
also the derivative $\dpr P/\dpr c$ vanishes at $c^*$.

Recall that, given two polynomials
\be\label{a.2}
P_{1}(c)=\sum_{i=0}^{n}a_{n-i}c^{i}, \qquad
P_{2}(c)=\sum_{i=0}^{m}b_{m-i}c^{i},
\ee
with $n,m\geq1$, the \emph{Sylvester matrix of} $P_{1},P_{2}$ is an
$n+m$ square matrix where the columns $1$ to $m$
are formed by ``shifted sequences'' of the coefficients of $P_{1}$,
while the columns $m+1$ to $m+n$ are formed by ``shifted sequences''
of the coefficients of $P_{2}$, \ie
\be\label{a.3}
\Syl(P_{1},P_{2}):=\begin{pmatrix}
a_{0} & 0 & \ldots & 0 & b_{0} & 0 & \ldots & 0 \cr
a_{1} & a_{0} & \ldots & 0 & b_{1} & b_{0} & \ldots & 0 \cr
\vdots & \vdots & \ddots & \vdots & \vdots & \vdots & \ddots & \vdots \cr
0 & 0 & \ldots & a_{n-1} & 0 & 0 & \ldots & b_{m-1} \cr
0 & 0 & \ldots & a_n & 0 & 0 & \ldots & b_{m}
\end{pmatrix},
\ee
and the \emph{resultant $R(P_{1},P_{2})$ of} $P_{1},P_{2}$
is defined as the determinant of the Sylvester matrix.

%%%%%%%%%%%%%%%%%%%%%%%%%%%%%%%%%%%%%%%%%%%%%%%%%%%%%%%%%%%%%%%%%%%%%%%%%
\begin{lemma}\label{lem:18}
Let $c_{1,1},\ldots,c_{1,n}$ and $c_{2,1},\ldots,c_{2,m}$ be the
complex roots of $P_{1},P_{2}$ respectively. Then 
\be\label{a.4}
R(P_{1},P_{2})=a_{0}^m b_{0}^n\prod_{i=1}^n\prod_{j=1}^m(c_{1,i}-c_{2,j}).
\ee
\end{lemma}
%%%%%%%%%%%%%%%%%%%%%%%%%%%%%%%%%%%%%%%%%%%%%%%%%%%%%%%%%%%%%%%%%%%%%%%%%

A complete proof is performed for instance in \cite{VDW}.
In particular, Lemma \ref{lem:18} implies that two polynomials
have a common root if and only if $R(P_{1},P_{2})=0$.

Recall also that given a polynomial $P=P(c)$, the
\emph{discriminant $D(P)$ of} $P$ is the resultant of $P$
and its first derivative with respect to $c$, \ie
$D(P):=R(P,P')$, where $P':=\de P/\de c$.
Thus, a polynomial $P=P(a,c)$ of the form (\ref{a.1})
has a multiple root if and only if its discriminant is equal to zero.

Now let us consider the set
\be\label{a.5}
V:=\{a=(a_{0},\ldots,a_n)\in\RRR^{n+1}\,:\,P(a,c)
\mbox{ has a multiple root}\}.
\ee
The discriminant of $P(a,c)$ is a polynomial
in the parameters $a=(a_{0},\ldots,a_n)$ \ie
$D_P(a)=D(P)\in\RRR[a_{0},\ldots,a_n]$, hence we can write
\be\label{a.6}
V=\{a=(a_{0},\ldots,a_n)\in\RRR^{n+1}\,:\, D_P(a)=0\}.
\ee
Such a set is, by definition, a proper Zariski-closed
subset of $\RRR^{n+1}$.

As the complement of a proper Zariski-closed subset of
$\RRR^{n+1}$ is open and dense also in the Euclidean topology,
then Hypothesis \ref{hyp3} is generic.\qed
%%%%%%%%%%%%%%%%%%%%%%%%%%%%%%%%%%%%%%%%%%%%%%%%%%%%%%%%%%%%%%%%%%%%%%%%%

%%%%%%%%%%%%%%%%%%%%%%%%%%%%%%%%%%%%%%%%%%%%%%%%%%%%%%%%%%%%%%%%%%%%%%%%%
%%%%%%%%%%%%%%%%%%%%%%%%%%%%%%%%%%%%%%%%%%%%%%%%%%%%%%%%%%%%%%%%%%%%%%%%%
\zerarcounters
\section{Proof of Lemma \ref{lem:17}}
\label{app:B}
%%%%%%%%%%%%%%%%%%%%%%%%%%%%%%%%%%%%%%%%%%%%%%%%%%%%%%%%%%%%%%%%%%%%%%%%%
%%%%%%%%%%%%%%%%%%%%%%%%%%%%%%%%%%%%%%%%%%%%%%%%%%%%%%%%%%%%%%%%%%%%%%%%%

\noindent First we shall prove by induction on $k$
that for all $\th\in\Th_{k,0,\b_{0}}$ one has
\be\label{b.1}
|L(\th)|\leq M(k-\hhh_{i_{0}})-\left(1+\frac{\sss_{i_{0}}}{\qqq}\right),
\ee
for all $k\geq\hhh_{i_{0}}+1$.

For $k=\hhh_{i_{0}}+1$ one has
\be\label{b.2}
{\hknu \b {\hhh_{i_{0}}+1} 0}=-\frac{1}{C}
\!\!\!\!\!\!\!\!\!\!\!\!\!\!\!\!
\!\!\!\!\!\!\!\!
\sum_{\substack{s_1,j\geq0\\ m_0+\ldots+m_{i_0}=j\\ s_{1}\ppp+
m_0\hhh_0+\ldots+m_{i_0}\hhh_{i_{0}}
=\sss_{i_{0}}+1}}
\!\!\!\!\!\!\!\!\!\!\!\!\!\!\!\!
\!\!\!\!\!\!\!\!
J(j,m_{0},\ldots,m_{i_{0}},m)\,
Q_{s_{1},j}c_0^{m_0}\ldots c_{i_{0}}^{m_{i_0}},
\ee
so that any tree $\th$ contributing to $\b_{0}^{[\hhh_{i_{0}}+1]}$ 
has $s_{1}+1$ nodes and $j$ leaves,
hence $|L(\th)|=s_{1}+1+j$.

Notice that $\qqq\leq\hhh_0\leq\hhh_1\leq\ldots\leq\hhh_{i_{0}}$,
hence one has
\be\label{b.3}
|L(\th)|=1+s_{1}+j \leq1+\frac{\sss_{i_{0}}+1}{\qqq}.
\ee

Moreover for $k=\hhh_{i_{0}}+1$ the r.h.s. in (\ref{b.1}) is
equal to $2+\sss_{i_{0}}/\qqq$,
so that the bound (\ref{b.1}) holds, because one has $\qqq\geq1$.
Assume now that the bound (\ref{b.1}) holds for
all $k'<k$ and let us show that then it holds also for $k$.

We call $M_{0}=M\hhh_{i_{0}}+1+\sss_{i_{0}}/\qqq$,
so that the inductive hypothesis can be written as
\be\label{b.4}
|L(\th')|\leq M k(\th')-M_{0},
\ee
for all $\th'\in\Th_{k',0,\b_{0}}$, $k'< k$.

By the inductive hypothesis, we have (cf. Figure \ref{fig:5})
\be\label{b.5}
|L(\th)|\leq1+s_{1}+s_{0}-s_{0}'M_{0}+M\sum_{i=1}^{s_{0}'}k(\th_{i}),
\ee
for suitable $\th_1,\ldots,\th_{s_0'}$ depending on $\th$.

Let us set $m:=k-\hhh_{i_{0}}\geq1$. Hence, \emph{via} the
conditions (\ref{eq:6.5}) we can write (\ref{b.5}) as
\be\label{b.6}
|L(\th)|\leq1+s_{1}+s_{0}-s_{0}'M_{0}+M\left(\sss_{i_{0}}+m-s_{1}\ppp-
\sum_{i=0}^{i_0}s_{0,i}\hhh_{i}\right).
\ee

Hence we shall prove that
\be\label{b.7}
1+s_{1}+s_{0}-s_{0}'M_{0}+M\left(\sss_{i_{0}}+m-s_{1}\ppp-
\sum_{i=0}^{i_0}s_{0,i}\hhh_{i}\right)
\leq mM-1-\frac{\sss_{i_{0}}}{\qqq},
\ee
or, in other words
\be\label{b.8}
\left(s_{1}\ppp+\sum_{i=0}^{i_0}s_{0,i}\hhh_{i}\right)M+s_{0}'M_0
\geq\sss_{i_{0}}
M+s_{0}+s_{1}+\frac{\sss_{i_{0}}}{\qqq}+2,
\ee
for all $s_{0},s_{0}',s_{1}\geq0$ admitted by conditions (\ref{eq:6.5}).

First of all for $s_{0}'=0$ by the first condition in (\ref{eq:6.5})
we have $s_{1}\ppp+s_{0,0}\hhh_0+\ldots+s_{0,i_0}\hhh_{i_{0}}=
\sss_{i_{0}}+m$.
Moreover $(s_{1}+s_{0})\qqq\leq s_{1}\ppp+s_{0,i_0}\hhh_0+\ldots+
s_{0,i_0}\hhh_{i_{0}}=\sss_{i_{0}}+m$,
hence
\be\label{b.9}
s_{1}+s_{0}\leq \frac{\sss_{i_{0}}+m}{\qqq},
\ee
so that one obtain (\ref{b.8}) if
\be\label{b.10}
m M\geq2 \frac{\sss_{i_{0}}}{\qqq}+2+\frac{m}{\qqq},
\ee
hence one needs
\be\label{b.11}
m\left(2 \frac{\sss_{i_{0}}}{\qqq}+3\right)\geq
2 \frac{\sss_{i_{0}}}{\qqq}+2+\frac{m}{\qqq},
\ee
that is satisfied for all $m\geq1$.

For $s_{0}'=1$ the first conditions (\ref{eq:6.5}) can be written as
$s_{1}\ppp+s_{0,0}\hhh_0+\ldots+(s_{0,i_0}+1)\hhh_{i_{0}}=
\sss_{i_{0}}+n$,
so that
\be\label{b.12}
s_{1}+s_{0}\leq \frac{\sss_{i_{0}}-\hhh_{i_{0}}+n}{\qqq}.
\ee

Hence we obtain (\ref{b.8}) if
\be\label{b.13}
nM\geq \frac{\sss_{i_{0}}+n-\hhh_{i_{0}}}{\qqq}+1,
\ee
and again (\ref{b.13}) is satisfied because $n=k-k_1\ge 1$.

Finally for $s_{0}'\geq2$ the first condition in (\ref{eq:6.5}) 
can be written
$s_{1}\ppp+s_{0,0}\hhh_0+\ldots+(s_{0,i_0}+s_{0}')
\hhh_{i_{0}}\geq\sss_{i_{0}}$,
so that
$s_{1}+s_{0}< s_{1}+s_{0}+s_{0}'\leq \frac{\sss_{i_{0}}}{\qqq}$,
and we obtain (\ref{b.8}) by requiring
\be\label{b.14}
\sss_{i_0} M+s_{0}'\left(\frac{\sss_{i_{0}}}{\qqq}+1\right)\geq\sss_{i_{0}}
M+2\frac{\sss_{i_{0}}}{\qqq}+2,
\ee
that is satisfied as we are assuming $s_{0}'\geq2$.

This exhausts the discussion
over all the choices of $s_{0},s_{0}',s_{1}$.

Let us show now that
\be\label{b.15}
|L(\th)|\leq M k-1,
\ee
for all $\th\in\Th_{k,\n,f}$, $f=\tilde{\b},B$, $k\geq\ppp$.

Again recall that a tree $\th\in\Th_{k,\n,f}$ contributes
to $f_{\n}^{[k]}$ with $f=\tilde{\b},B$, so that
the bound (\ref{b.15}) is trivially satisfied for $k=\ppp$
because one has $|L(\th)|=1$.

Let us suppose now that the bound holds for all $\ppp< k'<k$;
again we shall prove that then it holds also for $k$.

Recall that a tree contributing to ${\hknu f k \n}$ is of the
form depicted in Figure \ref{fig:6}, where $s_{0,\aaa}$
is the number of the lines exiting a leaf with leaf label $\aaa$
and entering $\vvv_{0}$, $s_0=s_{0,0}+\ldots+s_{0,i_0}$,
 $s_{1}$ is the number of the lines
exiting a node and entering $\vvv_{0}$,
and $s_{0}',s_{1}'$ are the graph elements entering $\vvv_{0}$
with component label $\b_{0}$ and $f$ respectively.
Hence, by the inductive hypothesis and by the bound (\ref{b.1}),
we have
\be\label{b.16}
|L(\th)|\leq1+s_{0}+s_{1}-s_{0}'M_{0}-s_{1}'+M
\sum_{i=1}^{s_{0}'+s_{1}'}k(\th_{i}),
\ee
for suitable $\th_1,\ldots,\th_{s_0'}$ depending on $\th$.

Let us supposte first $b_{\vvv_{0}}=1$; thus,
\emph{via} the first condition in (\ref{eq:6.2}), we
have to prove the bound
\be\label{b.17}
1+s_{0}+s_{1}+M(k-\ppp-s_{0,0}\hhh_{0}-\ldots-s_{0,i_{0}}
\hhh_{i_{0}}-s_{1}\ppp)-s_{0}'M_{0}-s_{1}'\leq M k-1,
\ee
or, in other words,
\be\label{b.18}
\sum_{i=0}^{i_{0}}s_{0,i}(M\hhh_{i}-1)+s_{1}(M\ppp-1)+M
\ppp+s_{0}'M_{0}+s_{1}'\geq2,
\ee
and this is obviously satisfied as $M\hhh_{i},M\ppp\geq3$.

Finally if $b_{\vvv_{0}}=0$ we have

\be\label{b.19}
\sum_{i=1}^{s_{1}'}k(\th_{i})=k-s_{1}\ppp,\qquad
s_{0,0}+\ldots+s_{0,i_{0}}+s_{0}'=0,\qquad s_{1}+s_{1}'\geq2,
\ee
so that, by the second condition in
(\ref{eq:6.2}), we have to prove the bound
\be\label{b.20}
1+s_{1}+M(k-s_{1}\ppp)-s_{1}'\leq M k-1,
\ee
or, in other words $s_{1}(M\ppp-1)+s_{1}'\geq2$,
and again this is obviously satisfied as $s_{1}+s_{1}'\geq2$
and $M\ppp>1$.

%%%%%%%%%%%%%%%%%%%%%%%%%%%%%%%%%%%%%%%%%%%%%%%%%%%%%%%%%%%%%%%%%%%%%%%%%%
%%%%%%%%%%%%%%%%%%%%%%%%%%%%%%%%%%%%%%%%%%%%%%%%%%%%%%%%%%%%%%%%%%%%%%%%%%
% References
%%%%%%%%%%%%%%%%%%%%%%%%%%%%%%%%%%%%%%%%%%%%%%%%%%%%%%%%%%%%%%%%%%%%%%%%%%
%%%%%%%%%%%%%%%%%%%%%%%%%%%%%%%%%%%%%%%%%%%%%%%%%%%%%%%%%%%%%%%%%%%%%%%%%%

\end{document}